\documentclass[12pt,draft]{article}

\usepackage{amsthm,amssymb,latexsym,amsmath}

\newtheorem{theorem}{Theorem}[section]
\newtheorem{definition}{Definition}[section]
\newtheorem{lemma}{Lemma}[section]

\newtheorem*{remark}{{\it Remark}}

\newcommand{\nc}{\newcommand}
\nc{\C}{{\mathbb C}}
\nc{\R}{{\mathbb R}}
\nc{\HH}{{\mathbb H}}
\nc{\Z}{{\mathbb Z}}
\nc{\N}{{\mathbb N}}
\nc{\dd}{{\rm d}}
\nc{\DD}{{\rm D}}

\begin{document}

\title{Moduli spaces of self-dual connections over asymptotically locally
flat gravitational instantons}

\author{G\'abor Etesi\\
\small{{\it Department of Geometry, Mathematical Institute, Faculty of
Science,}}\\
\small{{\it Budapest University of Technology and Economics,}}\\
\small{{\it Egry J. u. 1, H \'ep., H-1111 Budapest, Hungary
\footnote{{\tt etesi@math.bme.hu} and {\tt etesi@ime.unicamp.br}}}}\\
Marcos Jardim\\
\small{{\it Instituto de Matem\'atica, Estat\'\i stica e
Computa\c{c}\~ao Cient\'\i fica,}}\\
\small{{\it Universidade Estadual de Campinas,}}\\
\small{{\it C.P. 6065, 13083-859, Campinas, SP, Brazil\footnote{{\tt
jardim@ime.unicamp.br}}}}}

\maketitle

\pagestyle{myheadings}
\markright{G.Etesi, M. Jardim: Yang--Mills instantons on ALF gravitational
instantons}

\thispagestyle{empty}

\begin{abstract}
We investigate Yang--Mills instanton theory over four dimensional
asymptotically locally flat (ALF) geometries, including gravitational
instantons of this type, by exploiting the existence of a natural smooth
compactification of these spaces introduced by Hausel--Hunsicker--Mazzeo.

First referring to the codimension 2 singularity removal theorem of
Sibner--Sibner and R\aa de we prove that given a smooth, finite energy, 
self-dual SU$(2)$ connection over a complete ALF space, its energy is 
congruent to a Chern--Simons invariant of the boundary three-manifold if 
the connection satisfies a certain holonomy condition at infinity and its 
curvature decays rapidly.

Then we introduce framed moduli spaces of self-dual connections over
Ricci flat ALF spaces. We prove that the moduli space of smooth,
irreducible, rapidly decaying self-dual connections obeying the 
holonomy condition with fixed finite energy and prescribed asymptotic 
behaviour on a fixed bundle is a finite dimensional manifold. We 
calculate its dimension by a variant of the Gromov--Lawson 
relative index theorem.

As an application, we study Yang--Mills instantons over the flat
$\R^3\times S^1$, the multi-Taub--NUT family, and the Riemannian
Schwarzschild space.
\end{abstract}

\centerline{\small{AMS Classification: Primary: 53C26; Secondary:
53C29, 58J20, 58J28, 83C15}}

\section{Introduction}

By a {\it gravitational instanton} one usually means a
connected, four dimensional complete hyper-K\"ahler Riemannian manifold.
In particular, these spaces have SU$(2)\cong {\rm Sp}(1)$ holonomy;
consequently, they are Ricci flat, and hence solutions of the Riemannian
Einstein's vacuum equation. 

The only compact, four dimensional hyper-K\"ahler spaces are, up to 
(universal) covering, diffeomorphic to the flat torus $T^4$ or a $K3$ 
surface. 

The next natural step would be to understand non-compact gravitational 
instantons. Compactness in this case should be replaced by the 
condition that the metric be complete and decay to the flat metric at 
infinity somehow such that the Pontryagin number of the manifold be finite.

Such open hyper-K\"ahler examples can be constructed as follows. Consider
a connected, orientable compact four-manifold $\overline{M}$ with
connected boundary $\partial\overline{M}$ which is a smooth
three-manifold. Then the open manifold
$M:=\overline{M}\setminus\partial\overline{M}$ has a
decomposition $M=K\cup W$ where $K$ is a compact subset and
$W\cong\partial\overline{M}\times\R^+$ is
an open annulus or neck. Parameterize the half-line $\R^+$ by $r$.
Assume $\partial\overline{M}$ is fibered over a base 
manifold $B$ with fibers $F$ and the complete hyper-K\"ahler metric $g$ 
asymptotically and locally looks like 
$g\sim\dd r^2+r^2g_B +g_F$. In other words the base $B$ of the fibration 
blows up locally in a Euclidean way as $r\rightarrow\infty$, while the volume 
of the fiber remains finite. By the curvature decay, $g_F$ must be flat, 
hence $F$ is a connected, compact, orientable, flat manifold. On induction 
of the dimension of $F$, we can introduce several cases of increasing 
transcendentality, using the terminology of Cherkis and Kapustin \cite{che}:
\begin{itemize}
\item[(i)] $(M,g)$ is ALE (asymptotically locally Euclidean)
if $\dim F=0$;

\item[(ii)] $(M,g)$ is ALF (asymptotically locally flat) if $\dim F=1$,
in this case necessarily $F\cong S^1$ must hold;

\item[(iii)] $(M,g)$ is ALG (this abbreviation by induction) if $\dim
F=2$, in this case $F\cong T^2$;

\item[(iv)] $(M,g)$ is ALH if $\dim F=3$, in this case $F$ is
diffeomorphic to one of the six flat orientable three-manifolds.

\end{itemize}
Due to their relevance in quantum gravity or recently rather in
low-energy supersymmetric solutions of string theory and, last but
not least, their mathematical beauty, there has been some effort to
classify these spaces over the past decades.
Trivial examples for any class is provided by the space 
$\R^{4-\dim F}\times F$ with its flat product metric.

The first two non-trivial, infinite families were discovered by 
Gibbons and Hawking in 1976 \cite{gib-haw} in a rather explicit form. One 
of these families are the $A_k$ ALE or multi-Eguchi--Hanson spaces. In 
1989, Kronheimer gave a full classification
of ALE spaces \cite{kro} constructing them as minimal resolutions of
$\C^2/\Gamma $ where $\Gamma\subset{\rm SU}(2)$ is a finite subgroup i.e.,
$\Gamma$ is either a cyclic group $A_k$, $k\geq 0$, dihedral group $D_k$
with $k>0$, or one of the exceptional groups $E_l$ with $l=6,7,8$.

The other infinite family of Gibbons and Hawking is the $A_k$ ALF or
multi-Taub--NUT family. Recently another $D_k$ ALF family has been
constructed by Cherkis and Kapustin \cite{che-kap2} and in a more 
explicit form by Cherkis and Hitchin \cite{che-hit}. 

Motivated by string theoretical considerations, Cherkis and
Kapustin have suggested a classification scheme for ALF spaces as
well as for ALG and ALH \cite{che} although they relax the above 
asymptotical behaviour of the metric in these later two cases in order to 
obtain a natural classification. They claim that the $A_k$ and
$D_k$ families with $k\geq 0$ exhaust the ALF geometry (in this
enumeration $D_0$ is the Atiyah--Hitchin manifold). For the ALG case if 
we suppose that these spaces arise by deformations of elliptic fibrations 
with only one singular fiber, it is conjectured that the possibilities are 
$D_k$ with $0\leq k\leq 5$ (cf. \cite{che-kap1}) and $E_l$ with $l=6,7,8$. 
An example for a non-trivial ALH space is the minimal resolution of 
$(\R\times T^3)/\Z_2$. The trouble 
is that these spaces are more transcendental as $\dim F$ increases hence 
their constructions, involving twistor theory, Nahm transform, etc. are less 
straightforward and explicit.

To conclude this brief survey, we remark that the restrictive 
hyper-K\"ahler assumption on the metric, which appeared to be relevant in 
the more recent string theoretical investigations, excludes some examples 
considered as gravitational instantons in the early eighties. An important 
non-compact example which satisfies the ALF condition is for instance the 
Riemannian Schwarzschild space, which is Ricci flat but not 
hyper-K\"ahler \cite{haw}. For a more complete list of such ``old'' 
examples cf. \cite{egu-gil-han}.

From Donaldson theory we learned that the moduli spaces of SU$(2)$ instantons
over compact four-manifolds encompass lot of information about the original
manifold hence understanding SU$(2)$ instantons over gravitational instantons
also might be helpful in their classification. On the compact examples 
$T^4$ and the $K3$'s, Yang--Mills instantons can be studied via 
the usual methods, especially the celebrated Hitchin--Kobayashi 
correspondence. The full construction of SU$(2)$
instantons in the hyper-K\"ahler ALE case was carried out in important papers by
Nakajima \cite{nak} and Kronheimer--Nakajima in 1990 \cite{kro-nak}. However,
the knowledge regarding moduli spaces of instantons over non-trivial ALF 
spaces is rather limited, even in the hyper-K\"ahler ALF case, due to 
analytical difficulties. One has only sporadic examples of explicit 
solutions (cf. e.g. \cite{ete-hau2,ete-hau3}). Also very little is known 
about instanton theory over the Riemannian Schwarzschild space 
\cite{cha-duf,ete-hau1}. The only well studied case is the flat 
$\R^3\times S^1$ space; instantons over this space, also known as 
{\it calorons}, have been extensively studied in the literature, cf. 
\cite{bru-baa,bru-nog-baa,nye}. Close to nothing is known about instantons 
over non-trivial ALG and ALH geometries.

Studying Yang--Mills instanton moduli spaces over ALF spaces is certainly
interesting not only because understanding the reducible solutions already
leads to an encouraging topological classification result in the 
hyper-K\"ahler case \cite{ete}, but also due to their physical 
significance. In this paper we set the foundations for a general theory of 
Yang--Mills instantons over ALF spaces in the broad sense adopted in the 
eighties i.e., including not hyper-K\"ahler examples, too.

In Section \ref{sectwo} we exploit the existence of a natural smooth
compactification $X$ of an ALF space introduced by
Hausel--Hunsicker--Mazzeo \cite{hau-hun-maz}. Working over this
compact space, the asymptotical behaviour of any finite energy
connection over an ALF space can be analyzed by the codimension 2
singularity removal theorem of Sibner--Sibner \cite{sib-sib} and
R\aa de \cite{rad}. This guarantees the existence of a locally flat
connection $\nabla_\Gamma$ with fixed constant holonomy in infinity
to which the finite energy connection converges. First we prove in
Section \ref{sectwo} that the energy of a smooth, self-dual SU$(2)$ 
connection of finite energy which satisfies a certain holonomy condition 
(cf. condition (\ref{holonomia}) here) and has rapid curvature decay (in 
the sense of condition (\ref{lecsenges}) in the paper), is congruent to 
the Chern--Simons invariant $\tau_N(\Gamma_\infty )$ of the boundary $N$ of 
the ALF space (Theorem \ref{spektrum1} here). If the holonomy condition 
holds then $\nabla_\Gamma$ is in fact flat and $\Gamma_\infty$ is a fixed 
smooth gauge for the limiting flat connection restricted to the boundary. 
The relevant holonomy condition can be replaced by a simple topological 
criterion on the infinity of the ALF space, leading to a more explicit 
form of this theorem (cf. Theorem \ref{spektrum2}).

Then in Section \ref{secthree} we introduce framed instanton moduli
spaces $\mathcal{M}(e,\Gamma )$ of smooth, irreducible, rapidly decaying 
self-dual SU$(2)$ connections, obeying the holonomy condition, with fixed 
energy $e<\infty$ and asymptotical behaviour described by the flat
connection $\nabla_\Gamma$ on a fixed bundle. Referring to a variant
of the Gromov--Lawson relative index theorem \cite{gro-law} (cf.
Theorem \ref{index} here) we will be able to demonstrate that a
framed moduli space over a Ricci flat ALF space is either empty or forms a 
smooth manifold of dimension
\[\dim\mathcal{M}(e,\Gamma ) =8\left( e+\tau_N(\Theta_\infty )-\tau_N
(\Gamma_\infty )\right) -3b^-(X)\] 
where $\Theta_\infty$ is the restriction to $N$ of the trivial flat 
connection $\nabla_\Theta$ in some smooth gauge and $b^-(X)$ is the rank of 
the negative definite part of the intersection form of the 
Hausel--Hunsicker--Mazzeo compactification (cf. Theorem \ref{moduluster}).

In Section \ref{secfour} we apply our results on three classical
examples, obtaining several novel facts regarding instantons over
them.

First, we prove in Theorem \ref{kaloron} that any smooth, finite
energy caloron over $\R^3\times S^1$ automatically satisfies our holonomy
condition, has integer energy $e\in\N$ if it decays rapidly and that the 
dimension of the moduli space in this case is $8e$, in agreement with 
\cite{bru-baa}. These moduli spaces are non-empty for all positive integer $e$ 
\cite{bru-nog-baa}.

For the canonically oriented multi-Taub--NUT spaces, we show that
the dimension of the framed moduli of smooth, irreducible, rapidly 
decaying anti-self-dual connections satisfying the holonomy condition is
divisible by 8 (cf. Theorem \ref{taub}).
%This contradicts the expectations that the moduli space should have a
%four dimensional component which is a copy of the original manifold as in the ALE case.
Known explicit solutions \cite{ete-hau3} show that at least a
few of these moduli spaces are actually non-empty.

Finally, we consider the Riemannian Schwarzschild case, and prove in
Theorem \ref{schwarz} that all smooth finite energy instantons obey
the holonomy condition, have integer energy $e$ if they decay rapidly and 
the dimension is $8e-3$. Moreover, this moduli space is surely non-empty at 
least for $e=1$. We also enumerate the remarkably few known explicit 
solutions \cite{cha-duf,ete-hau1}, and observe that these admit 
deformations.

Section \ref{secfive} is an Appendix containing the proof of the
relative index theorem used in the paper, Theorem  \ref{index}.

%%%%%%%%%%%%%%%%%%%%%%%%%%%%%%%%%%%%%%%%%%%%%%%%%%%%%%%%%%%%%%%%%%%%%%%
%%%%%%%%%%%%%%%%%%%%%%%%%%%%%%%%%%%%%%%%%%%%%%%%%%%%%%%%%%%%%%%%%%%%%%%

\section{The spectrum of the Yang--Mills functional}\label{sectwo}

In this Section we prove that the spectrum of the Yang--Mills
functional evaluated on self-dual connections satisfying a certain analytical 
and a topological condition over a complete ALF manifold is ``quantized'' 
by the Chern--Simons invariants of the boundary. First, let us carefully 
define the notion of ALF space used in this paper, and describe its useful
topological compactification, first used in \cite{hau-hun-maz}.

Let $(M,g)$ be a connected, oriented Riemannian four-manifold. This
space is called an {\it asymptotically locally flat (ALF) space} if
the following holds. There is a compact subset $K\subset M$ such
that $M\setminus K =W$ and $W\cong N\times\R^+$, with $N$ being a
connected, compact, oriented three-manifold without boundary
admitting a smooth $S^1$-fibration
\begin{equation}
\pi :N\stackrel{F}{\longrightarrow}B_\infty
\label{perem}
\end{equation}
whose base space is a compact Riemann surface $B_\infty$.
For the smooth, complete Riemannian metric $g$ there exists a
diffeomorphism $\phi :N\times\R^+\rightarrow W$ such that
\begin{equation}
\phi^*(g\vert_W)=\dd r^2+r^2(\pi^*g_{B_\infty})'+h'_F
\label{aszimptotika}
\end{equation}
where $g_{B_\infty}$ is a smooth metric on $B_\infty$, $h_F$ is a
symmetric 2-tensor on $N$ which restricts to a metric along the
fibers $F\cong S^1$ and $(\pi^*g_{B_\infty})'$ as well as $h'_F$ are
some finite, bounded, smooth extensions of $\pi^*g_{B_\infty}$ and $h_F$ 
over $W$, respectively. That is, we require $(\pi^*g_{B_\infty})'(r)\sim
O(1)$ and $h'_F(r)\sim O(1)$ and the extensions for $r<\infty$ preserve 
the properties of the original fields. Furthermore, 
we impose that the curvature $R_g$ of $g$ decays like
\begin{equation}
\vert\phi^*(R_g\vert_W)\vert\sim O(r^{-3}).
\label{gorbulet}
\end{equation}
Here $R_g$ is regarded as a map
$R_g:C^\infty (\Lambda^2M)\rightarrow C^\infty (\Lambda^2M)$
and its pointwise norm is calculated accordingly in an orthonormal frame. 
Hence the Pontryagin number of our ALF spaces is finite.

Examples of such metrics are the natural metric on $\R^3\times S^1$
which is in particular flat; the multi-Taub--NUT family
\cite{gib-haw} which is moreover hyper-K\"ahler or the Riemannian
Schwarzschild space \cite{haw} which is in addition Ricci flat only.
For a more detailed description of these spaces, cf. Section
\ref{secfour}.

We construct the compactification $X$ of $M$ simply by shrinking all
fibers of $N$ into points as $r\to\infty$ like in \cite{hau-hun-maz}. We 
put an orientation onto $X$ induced by the orientation of the original 
$M$. The space $X$ is then a connected, oriented, smooth four-manifold
without boundary. One clearly obtains a decomposition $X\cong M\cup
B_\infty$, and consequently we can think of $B_\infty$ as a smoothly
embedded submanifold of $X$ of codimension 2. For example, for $\R^3\times 
S^1$ one finds $X\cong S^4$ and $B_\infty$ is an embedded $S^2$ 
\cite{ete}; for the multi-Taub--NUT space with the orientation induced by 
one of the complex structures, $X$ is the connected sum of $s$ copies of
$\overline{\C P}^2$'s ($s$ refers to the number of NUTs) and
$B_\infty$ is homeomorphic to $S^2$ providing a generator of the
second cohomology of $X$; in case of the Riemannian Schwarzschild
geometry, $X\cong S^2\times S^2$ and $B_\infty$ is again $S^2$, also
providing a generator for the second cohomology (cf. \cite{ete,hau-hun-maz}).

Let $M_R:=M\setminus (N\times (R,\infty ))$ be the truncated
manifold with boundary $\partial M_R\cong N\times\{ R\}$. Taking
into account that $W\cong N\times (R,\infty )$, a normal
neighbourhood $V_R$ of $B_\infty$ in $X$ has a model like $V_R\cong
N\times (R,\infty ]/\sim$ where $\sim$ means that $N\times\{\infty\}$ 
is pinched into $B_\infty$. We obtain $W =V_R\setminus V_\infty$, with
$V_\infty\cong B_\infty$. By introducing the parameter 
$\varepsilon :=R^{-1}$ we have another model $V_\varepsilon$ provided by 
the fibration 
\[V_\varepsilon\stackrel{B_\varepsilon^2}{\longrightarrow }B_\infty\]
whose fibers are two-balls of radius $\varepsilon$. In this second
picture we have the identification $V_0\cong B_\infty$, so that the
end $W$ looks like
\begin{equation}
V^*_\varepsilon :=V_\varepsilon\setminus V_0.
\label{kornyezet}
\end{equation}
Choosing a local coordinate patch $U\subset B_\infty$, then locally
$V_\varepsilon\vert_U\cong U\times B^2_\varepsilon$ and
$V^*_\varepsilon\vert_U\cong U\times (B^2_\varepsilon\setminus\{ 0\} )$.
We introduce local coordinates $(u,v)$ on $U$ and polar coordinates
$(\rho ,\tau )$ along the discs $B^2_\varepsilon$ with
$0\leq \rho <\varepsilon$ and $0\leq\tau <2\pi $. Note that in fact 
$\rho =r^{-1}$ is a global coordinate over the whole $V_\varepsilon\cong 
V_R$. For simplicity we denote $V_\varepsilon\vert_U$ as $U_\varepsilon$ and 
will call the set
\begin{equation}
U^*_\varepsilon :=U_\varepsilon\setminus U_0
\label{elemikornyezet}
\end{equation}
an {\it elementary neighbourhood}. Clearly, their union covers the
end $W$. In this $\varepsilon$-picture we will use the notation
$\partial M_\varepsilon\cong N\times\{\varepsilon\}$ for the
boundary of the truncated manifold $M_\varepsilon =M\setminus
(N\times (0,\varepsilon ))$ and by a slight abuse of notation we
will also think of the end sometimes as $W\cong\partial
M_\varepsilon\times (0, \varepsilon )$.

We do not expect the complete ALF metric $g$ to extend over this
compactification, even conformally. However the ALF property
(\ref{aszimptotika}) implies that we can suppose the existence of a smooth
positive function $f\sim O(r^{-2})$ on $M$ such that the rescaled metric
$\tilde{g}:=f^2g$ extends smoothly as a tensor field over $X$ (i.e., a
smooth Riemannian metric degenerated along the singularity set
$B_\infty$). In the vicinity of the singularity we find
$\tilde{g}\vert_{V_\varepsilon}=\dd\rho^2+\rho^2
(\pi^*g_{B_\infty})''+\rho^4h''_F$ via (\ref{aszimptotika}) consequently
we can choose the coordinate system $(u,v,\rho ,\tau)$ on
$U_\varepsilon$ such that $\{\dd u,\dd v, \dd\rho ,\dd\tau\}$ 
forms an oriented frame on $T^*U^*_\varepsilon$ and with some bounded, 
finite function $\varphi$, the metric looks like
$\tilde{g}\vert_{U^*_\varepsilon}=\dd\rho^2+\rho^2\varphi (u,v,\rho )(\dd 
u^2+\dd v^2)+\rho^4(\dd\tau^2+\dots)$. Consequently we find that
\begin{equation}
{\rm Vol}_{\tilde{g}}(V_\varepsilon)\sim O(\varepsilon^5).
\label{terfogat}
\end{equation}

We will also need a smooth regularization of $\tilde{g}$. Taking a 
monotonously increasing smooth function $f_\varepsilon$ supported in 
$V_{2\varepsilon}$ and equal to 1 on $V_\varepsilon$ such that $\vert\dd 
f_\varepsilon\vert\sim O(\varepsilon^{-1})$, as well as picking up 
a smooth metric $h$ on $X$, we can regularize $\tilde{g}$ by introducing the 
smooth metric
\begin{equation}
\tilde{g}_\varepsilon :=(1-f_\varepsilon )\tilde{g} +f_\varepsilon h
\label{metrika}
\end{equation}
over $X$. It is clear that $\tilde{g}_0$ and $\tilde{g}$ agree on $M$.

Let $F$ be an SU$(2)$ vector bundle over $X$ endowed with a fixed 
connection $\nabla_\Gamma$ and an invariant fiberwise scalar 
product. Using the rescaled-degenerated metric $\tilde{g}$, define 
Sobolev spaces $L^p_{j, \Gamma}(\Lambda^*X\otimes F)$ with $1<p\leq\infty$ 
and $j=0,1,\dots$ as the completion of $C_0^\infty (\Lambda^*X\otimes F)$, 
smooth sections compactly supported in $M\subset X$, with respect to the norm
\begin{equation}
\Vert\omega\Vert_{L^p_{j ,\Gamma}(X)}:=
\left(\lim\limits_{\varepsilon\rightarrow
0}\sum\limits_{k=0}^j\Vert\nabla_\Gamma^k 
\omega\Vert^p_{L^p(M_\varepsilon,\tilde{g}\vert_{M_\varepsilon})}\right)
^{\frac{1}{p}}
\label{normak}
\end{equation}
where
\[\Vert\nabla^k_\Gamma\omega\Vert^p_{L^p(M_\varepsilon,\tilde{g}
\vert_{M_\varepsilon})}=
\int\limits_{M_\varepsilon}\vert\nabla^k_\Gamma\omega\vert^p*_{\tilde{g}}1.\]
Throughout this paper, Sobolev norms of this kind will be used
unless otherwise stated. We will write simply $L^p$ for 
$L^p_{0,\Gamma}$. Notice that every 2-form with finite $L^2$-norm over 
$(M,g)$ will also belong to this Sobolev space, by conformal invariance 
and completeness.

Next, we collect some useful facts regarding the Chern--Simons
functional. Let $E$ be an smooth SU$(2)$ bundle over $M$. Since 
topological $G$-bundles over an open four-manifold are classified by 
$H^2(M,\pi_1(G))$, note that $E$ is necessarily trivial. Put a 
smooth SU$(2)$ connection $\nabla_B$ onto $E$. Consider the boundary $\partial
M_\varepsilon$ of the truncated manifold. The restricted bundle
$E\vert_{\partial M_\varepsilon}$ is also trivial. Therefore any
restricted SU$(2)$ connection $\nabla_B\vert_{\partial
M_\varepsilon}:=\nabla_{B_\varepsilon}$ over $E\vert_{\partial
M_\varepsilon}$ can be identified with a smooth ${\mathfrak
s}{\mathfrak u}(2)$-valued 1-form $B_\varepsilon$. The {\it
Chern--Simons functional} is then defined to be
\[\tau_{\partial M_\varepsilon}(B_\varepsilon
):=-\frac{1}{8\pi^2}\int\limits_{\partial M_\varepsilon}{\rm
tr}\left(\dd B_\varepsilon\wedge
B_\varepsilon+\frac{2}{3}B_\varepsilon\wedge B_\varepsilon\wedge
B_\varepsilon\right) .\]
This expression is gauge invariant up to an integer. Moreover, the
representation space
\[\chi (\partial M_\varepsilon ):={\rm Hom}(\pi_1(\partial M_\varepsilon
), {\rm SU}(2))/{\rm SU}(2)\]
is called the {\it character variety} of $\partial M_\varepsilon\cong N$
and parameterizes the gauge equivalence classes of smooth flat SU$(2)$ 
connections over $N$.

\begin{lemma} Fix an $0<\rho <\varepsilon$ and let
$\nabla_{A_\rho}=\dd +A_\rho$ and
$\nabla_{B_\rho}=\dd + B_\rho$ be two smooth SU$(2)$
connections in a fixed smooth gauge on the trivial SU$(2)$ bundle
$E\vert_{\partial M_\rho}$. Then there is a constant 
$c_1=c_1(B_\rho )>0$, depending on $\rho$ only through $B_\rho$, 
such that
\begin{equation}
\vert \tau_{\partial M_\rho}(A_\rho )-\tau_{\partial
M_\rho}(B_\rho )\vert \leq c_1\Vert
A_\rho -B_\rho\Vert_{L^2_{1,B_\rho}(\partial M_\rho )}
\label{cs}
\end{equation}
that is, the Chern--Simons functional is continuous in the 
$L^2_{1,B_\rho}$ norm.

Moreover, for each $\rho$, $\tau_{\partial M_\rho}(A_\rho)$ is constant on the
path connected components of the character variety $\chi (\partial M_\rho )$.
\label{cslemma}
\end{lemma}

\noindent{\it Proof.} The first observation follows from the identity
\[\tau_{\partial M_\rho}(A_\rho )-\tau_{\partial
M_\rho}(B_\rho )=\]
\[-\frac{1}{8\pi^2}\int\limits_{\partial M_\rho}{\rm
tr}\left( (F_{A_\rho}+F_{B_\rho})\wedge (A_\rho -B_\rho 
)-\frac{1}{3}(A_\rho -B_\rho )\wedge (A_\rho -B_\rho )\wedge
(A_\rho -B_\rho )\right)\]
which implies that there is a constant $c_0=c_0(\rho ,B_\rho )$ such that
\[\vert \tau_{\partial M_\rho}(A_\rho )-\tau_{\partial
M_\rho}(B_\rho )\vert\leq c_0\Vert A_\rho 
-B_\rho\Vert_{L^{\frac{3}{2}}_{1,B_\rho}(\partial M_\rho )}\] 
that is, the Chern--Simons functional is continuous in the
$L^{\frac{3}{2}}_{1,B_\rho}$ norm. 
A standard application of H\"older's inequality on $(\partial
M_\rho ,\tilde{g}\vert_{\partial M_\rho})$ then yields
\[\Vert A_\rho
-B_\rho\Vert_{L^{\frac{3}{2}}_{1,B_\rho}(\partial M_\rho )}\leq 
\sqrt{2}\left({\rm Vol}_{\tilde{g}\vert_{\partial M_\rho}}(\partial 
M_\rho )\right)^{\frac{1}{6}}\Vert A_\rho
-B_\rho\Vert_{L^2_{1,B_\rho}(\partial M_\rho )}.\]

The metric locally looks 
like $\tilde{g}\vert_{\partial M_\rho\cap U^*_\varepsilon} 
=\rho^2\varphi (\dd u^2+\dd v^2)+\rho^4(\dd\tau^2+2h_{\tau ,u}\dd\tau\dd 
u+\dots )$ with $\varphi$ and $h_{\tau 
,u}$, etc. being bounded functions of $(u,v,\rho )$ and $(u,v,\rho ,\tau 
)$ respectively hence the metric coefficients as well as the volume of 
$(\partial M_\rho ,\tilde{g}\vert_{\partial M_\rho})$ are bounded 
functions of $\rho$ consequently we 
can suppose that $c_1$ does not depend explicitly on $\rho$. 

Concerning the second part, assume $\nabla_{A_\rho}$ and $\nabla_{B_\rho}$
are two smooth, flat connections belonging to the same path connected component 
of $\chi(\partial M_\rho)$. Then there is a continuous path
$\nabla_{A^t_\rho}$ with $t\in [0,1]$ of flat connections connecting the 
given flat connections. Out of this we construct a connection $\nabla_A$ 
on $\partial M_\rho\times [0,1]$ given by $A:=A^t_\rho +0\cdot\dd t$. 
Clearly, this connection is flat, i.e., $F_A=0$. The Chern--Simons theorem 
\cite{che-sim} implies that
\[\tau_{\partial M_\rho}(A_\rho )-\tau_{\partial
M_\rho}(B_\rho )=-\frac{1}{8\pi^2}\int\limits_{\partial
M_\rho\times [0,1]}{\rm tr} (F_A\wedge F_A)=0,\]
concluding the proof. $\Diamond$
\vspace{0.1in}

\noindent The last ingredient in our discussion is the fundamental theorem
of Sibner--Sibner \cite{sib-sib} and R\aa de \cite{rad} which allows us to
study the asymptotic behaviour of finite energy connections over an ALF
space. Consider a smooth (trivial) SU$(2)$ vector bundle $E$ over the ALF 
space $(M,g)$. Let $\nabla_A$ be a Sobolev connection on $E$ with 
finite energy, i.e. $F_A\in L^2 (\Lambda^2M\otimes{\rm End}E)$. 
Taking into account completeness of the ALF metric, we have 
$\vert F_A\vert (r)\to 0$ almost everywhere as $r\to\infty$. 
Thus we have a connection defined on $X$ away from a smooth, codimension 2 
submanifold $B_\infty\subset X$ and satisfies $\Vert 
F_A\Vert_{L^2(X)}<\infty$.

Consider a neighbourhood $B_\infty\subset V_\varepsilon$ and write
$V^*_\varepsilon$ to describe the end $W$ as in (\ref{kornyezet}).
Let $\nabla_\Gamma$ be an SU$(2)$ connection on 
$E\vert_{V^*_\varepsilon}$ which is locally flat and smooth. The 
restricted bundle $E\vert_{V^*_\varepsilon}$ is trivial, hence we can 
choose some global gauge such that $\nabla_A\vert_{V^*_\varepsilon}=\dd
+A_{V^*_\varepsilon}$ and $\nabla_\Gamma =\dd +\Gamma_{V^*_\varepsilon}$; 
we assume with some $j=0,1,\dots$ that $A_{V^*_\varepsilon}\in L^2_{j+1,\Gamma 
,loc}$. Taking into account 
that for the elementary neighbourhood $\pi_1(U^*_\varepsilon )\cong\Z$, 
generated by a $\tau$-circle, it is clear
that locally on $U^*_\varepsilon\subset {V^*_\varepsilon}$ we can choose a
more specific gauge $\nabla_\Gamma\vert_{U^*_\varepsilon}=\dd +\Gamma_m$
with a constant $m\in [0,1)$ such that \cite{sib-sib}
\begin{equation}
\Gamma_m=\begin{pmatrix}{\bf i}m & 0\cr
                      0 & -{\bf i}m
          \end{pmatrix}\dd\tau .
\label{lokalismertek}
\end{equation}
Here $m$ represents the {\it local holonomy} of the locally 
flat connection around the punctured discs of the space $U^*_\varepsilon$, 
see \cite{sib-sib}. It is invariant under gauge transformations modulo an 
integer. For later use we impose two conditions on this local holonomy. 
The embedding $i:U^*_\varepsilon\subset {V^*_\varepsilon}$ 
induces a group homomorphism $i_*: 
\pi_1(U^*_\varepsilon)\rightarrow\pi_1({V^*_\varepsilon})$.
It may happen that this homomorphism has non-trivial kernel. Let $l$ be
a loop in $U^*_\varepsilon$ such that $[l]$ generates 
$\pi_1(U^*_\varepsilon )\cong\Z$.

\begin{definition}
A locally flat connection $\nabla_{\Gamma}$ on $E\vert_{V^*_\varepsilon}$ 
is said to satisfy the {\rm  weak holonomy condition} if for all 
$U^*_\varepsilon\subset V^*_\varepsilon$ the restricted connection
$\nabla_{\Gamma}\vert_{U^*_\varepsilon}=\dd +\Gamma_m$ has trivial local 
holonomy whenever $l$ is contractible in ${V^*_\varepsilon}$, i.e.
\begin{equation}
[l]\in{\rm Ker}\: i_* ~~ \Longrightarrow ~~ m=0.
\label{holonomia}
\end{equation}
Additionally, $\nabla_{\Gamma}$ is said to satisfy the {\rm strong 
holonomy condition} if the local holonomy of any restriction 
$\nabla_\Gamma\vert_{U^*_\varepsilon}=\dd +\Gamma_m$ vanishes, i.e.
\begin{equation}
m=0.
\label{holonomia2}
\end{equation}
\label{whc}
\end{definition}

\noindent Clearly, $\nabla_\Gamma$ is {\it globally} a smooth, flat 
connection on $E\vert_{V^*_\varepsilon}$ if and only if the weak holonomy 
condition holds. Moreover, the strong holonomy condition implies the weak one. 
We are now in a position to recall the following fundamental 
regularity result \cite{rad,sib-sib}. 

\begin{theorem}{\rm (Sibner--Sibner, 1992 and R\aa de, 1994)} There
exist a constant $\varepsilon >0$ and a flat SU$(2)$
connection $\nabla_\Gamma\vert_{U^*_\varepsilon}$ on 
$E\vert_{U^*_\varepsilon}$ with a constant holonomy $m\in [0,1)$ 
such that on $E\vert_{U^*_\varepsilon}$ one can find a gauge
$\nabla_A\vert_{U^*_\varepsilon}=\dd +A_{U^*_\varepsilon}$ and
$\nabla_\Gamma\vert_{U^*_\varepsilon}=\dd +\Gamma_m$ with
$A_{U^*_\varepsilon}-\Gamma_m\in L^2_{1,\Gamma}(U^*_\varepsilon )$ such that 
the estimate
\[\Vert A_{U^*_\varepsilon}-\Gamma_m\Vert_{L^2_{1,\Gamma}(U^*_\varepsilon
)}\leq c_2\Vert F_A\Vert_{L^2(U_\varepsilon )}\]
holds with a constant $c_2=c_2(\tilde{g}\vert_{U_\varepsilon})>0$ 
depending only on the metric. $\Diamond$
\label{lokalissibner}
\end{theorem}
\noindent This theorem shows that any finite energy connection is always 
asymptotic to a flat connection at least locally. It is therefore convenient 
to say that the finite energy connection $\nabla_A$ {\it satisfies the weak} or
{\it the strong holonomy condition} if its associated asymptotic 
locally flat connection $\nabla_\Gamma$, in the sense of Theorem 
\ref{lokalissibner}, satisfies the corresponding condition in the sense 
of Definition \ref{whc}. We will be using this terminology.

This estimate can be globalized over the whole end
${V^*_\varepsilon}$ as follows.
Consider a finite covering $B_\infty =\cup_\alpha U_\alpha$ and denote the
corresponding punctured sets as $U^*_{\varepsilon, \alpha}\subset
{V^*_\varepsilon}$. These sets also give rise to a finite
covering of ${V^*_\varepsilon}$.
It is clear that the weak condition (\ref{holonomia}) is independent of
the index $\alpha$, since by Theorem \ref{lokalissibner}, $m$ is 
constant over all ${V^*_\varepsilon}$. Imposing (\ref{holonomia}), 
the local gauges $\Gamma_m$ on $U^*_{\varepsilon ,\alpha}$ extend smoothly 
over the whole $E\vert_{V^*_\varepsilon}$. That is, there is a smooth flat 
gauge $\nabla_\Gamma =\dd +\Gamma_{V^*_\varepsilon}$ over 
$E\vert_{V^*_\varepsilon}$
such that $\Gamma_{V^*_\varepsilon}\vert_{U^*_{\varepsilon 
,\alpha}}=\gamma^{-1}_\alpha
\Gamma_m\gamma_\alpha +\gamma^{-1}_\alpha\dd\gamma_\alpha$ with smooth 
gauge transformations $\gamma_\alpha : U^*_{\varepsilon 
,\alpha}\rightarrow {\rm SU}(2)$. This gauge is unique only up to an 
arbitrary smooth gauge transformation. Since this construction deals with 
the topology of the boundary (\ref{perem}) only, we can assume that these 
gauge transformations are independent of the (global) radial coordinate 
$0<\rho <\varepsilon$. Then we write this global gauge as
\begin{equation}
\nabla_A\vert_{V^*_\varepsilon}=\dd
+A_{V^*_\varepsilon},\:\:\:\:\:\nabla_\Gamma =\dd 
+\Gamma_{V^*_\varepsilon} .
\label{globalismertek}
\end{equation}
A comparison with the local gauges in Theorem \ref{lokalissibner} shows that
\begin{equation}
(A_{V^*_\varepsilon}-\Gamma_{V^*_\varepsilon} )\vert_{U^*_{\varepsilon ,
\alpha}}=\gamma^{-1}_\alpha (A_{U^*_{\varepsilon 
,\alpha}}-\Gamma_m)\gamma_\alpha .
\label{ujmerce}
\end{equation}
Applying Theorem \ref{lokalissibner} in all coordinate
patches and summing up over them we come up with 
\begin{equation}
\Vert A_{V^*_\varepsilon}-\Gamma_{V^*_\varepsilon}
\Vert_{L^2_{1,\Gamma}(V^*_\varepsilon 
)}\leq c_3\Vert F_A\Vert_{L^2(V_\varepsilon )},
\label{globalissibner}
\end{equation}
with some constant $c_3=c_3 (\tilde{g}\vert_{V_\varepsilon},\gamma_\alpha 
,\dd\gamma_\alpha )>0$. This is the globalized version of Theorem 
\ref{lokalissibner}. 

Let $A_\varepsilon$ and $\Gamma_\varepsilon$ be the restrictions of
the connection 1-forms in the global gauge (\ref{globalismertek}) to 
the boundary $E\vert_{\partial M_\varepsilon}$. The 
whole construction shows that $\Gamma_\varepsilon$ is
smooth and is independent of $\varepsilon$ consequently 
$\lim\limits_{\varepsilon\rightarrow 0}\Gamma_\varepsilon$
exists and is smooth. Note again that smoothness follows if and only if 
the weak holonomy condition (\ref{holonomia}) is satisfied. In this way
$\lim\limits_{\varepsilon\rightarrow 0}\tau_{\partial
M_\varepsilon}(\Gamma_\varepsilon)=\tau_N(\Gamma_0)$ also exists and gives 
rise to a {\it Chern--Simons invariant} of the boundary.

Assume that our finite energy connection $\nabla_A$ is smooth; then in the 
global gauge (\ref{globalismertek}) $A_{V^*_\varepsilon}$ is also smooth, 
hence $\tau_{\partial M_\varepsilon}(A_\varepsilon )$ also exists for 
$\varepsilon >0$. 

Next we analyze the behaviour $\tau_{\partial M_\varepsilon}(A_\varepsilon 
)$ as $\varepsilon$ tends to zero. First notice that the local 
flat gauge $\Gamma_m$ in (\ref{lokalismertek}) does not have radial component 
consequently $A_{V^*_\varepsilon}-\Gamma_{V^*_\varepsilon} =A_\varepsilon 
-\Gamma_\varepsilon +A_r$
where $A_r$ is the radial component of $A_{V^*_\varepsilon}$. Dividing the 
square of (\ref{cs}) by $\varepsilon >0$ and then integrating it we 
obtain, making use of (\ref{globalissibner}) that
\[\frac{1}{\varepsilon}\int\limits_{0}^\varepsilon\vert\tau_{\partial 
M_\rho}(A_\rho )-\tau_{\partial M_\rho}(\Gamma_\rho)\vert^2\dd\rho \leq 
\frac{c^2_1}{\varepsilon}\int\limits_0^\varepsilon\Vert A_\rho 
-\Gamma_\rho\Vert^2_{L^2_{1,\Gamma_\rho}(\partial M_\rho )}\dd\rho\leq\]
\[ \leq \frac{c^2_1}{\varepsilon}\int\limits_0^\varepsilon\left(\Vert 
A_\rho -\Gamma_\rho\Vert^2_{L^2_{1,\Gamma_\rho}(\partial M_\rho )}+\Vert 
A_r\Vert^2_{L^2_{1,\Gamma_\rho}(\partial M_\rho )}\right)\dd\rho 
\leq\frac{c^2_1}{\varepsilon}\:\Vert A_{V^*_\varepsilon}-
\Gamma_{V^*_\varepsilon}\Vert^2_{L^2_{1,\Gamma}(V^*_\varepsilon 
)}\leq \frac{(c_1c_3)^2}{\varepsilon}\Vert F_A\Vert^2_{L^2(V_\varepsilon )} .\]
Finite energy and completeness implies that $\Vert 
F_A\Vert_{L^2(V_\varepsilon )}$ vanishes as $\varepsilon$ 
tends to zero. However for our purposes we need a stronger decay assumption.
\begin{definition}
The finite energy SU$(2)$ connection $\nabla_A$ on the bundle $E$ over 
$M$ {\em decays rapidly} if its curvature satisfies
\begin{equation}
\lim\limits_{\varepsilon\rightarrow 0}\frac{\Vert
F_A\Vert_{L^2(V_\varepsilon )}}
{\sqrt{\varepsilon}}=\lim\limits_{R\rightarrow\infty}\sqrt{R}\:\Vert
F_A\Vert_{L^2(V_R,g\vert_{V_R})}=0
\label{lecsenges}
\end{equation}
along the end of the ALF space.
\end{definition}

\noindent Consequently for a rapidly decaying connection we obtain 
\[\lim\limits_{\varepsilon\rightarrow 
0}\frac{1}{\varepsilon}\int\limits_0^\varepsilon\vert\tau_{\partial
M_\rho}(A_\rho )-\tau_{\partial M_\rho}(\Gamma_\rho)\vert^2\dd\rho =0,\]
which is equivalent to
\begin{equation}
\lim\limits_{\varepsilon\rightarrow 0}\tau_{\partial
M_\varepsilon}(A_\varepsilon )=\tau_N(\Gamma_0).
\label{magnes}
\end{equation}

We are finally ready to state an energy identity for self-dual
connections. Let $\nabla_A$ be a smooth, self-dual, finite energy connection 
on the trivial SU$(2)$ bundle $E$ over an ALF space $(M,g)$:
\[F_A=*F_A,\:\:\:\:\:\Vert F_A\Vert^2_{L^2(M,g)}<\infty .\]
Assume it satisfies the weak holonomy condition (\ref{holonomia}).
In this case we can fix a gauge (\ref{globalismertek}) along
$V^*_\varepsilon$ and both $A_{V^*_\varepsilon}$ and 
$\Gamma_{V^*_\varepsilon}$ are smooth. Restrict $\nabla_A$ onto 
$E\vert_{M_\varepsilon}$ with
$\varepsilon >0$. Exploiting self-duality, an application of the
Chern--Simons theorem \cite{che-sim} along the boundary shows that
\[\Vert F_A\Vert^2_{L^2(M_\varepsilon
,g\vert_{M_\varepsilon})}\equiv\tau_{\partial
M_\varepsilon}(A_\varepsilon )\:\:\:\:\:\mbox{mod $\Z$}.\]
Moreover if the connection decays rapidly in the sense of 
(\ref{lecsenges}) then the right hand side has a limit (\ref{magnes}) 
therefore we have arrived at the following theorem:

\begin{theorem}
Let $(M,g)$ be an ALF space with an end $W\cong N\times\R^+$. Let $E$ be 
an SU$(2)$ vector bundle over $M$, necessarily trivial, with a 
smooth, finite energy, self-dual connection $\nabla_A$. If it satisfies 
the weak holonomy condition (\ref{holonomia}) and decays rapidly in the 
sense of (\ref{lecsenges}) then there exists a smooth flat 
SU$(2)$ connection $\nabla_\Gamma$ on $E\vert_W$ and a smooth flat gauge 
$\nabla_\Gamma =\dd +\Gamma_W$, unique up to a smooth gauge 
transformation, such that $\lim\limits_{r\rightarrow\infty}
\Gamma_W\vert_{N\times\{ r\}}=\Gamma_\infty$ exists, is smooth and 
\[\Vert F_A\Vert^2_{L^2(M,g)}\equiv\tau_N
(\Gamma_\infty )\:\:\:\:\:\mbox{{\rm mod} $\Z$}.\]
That is, the energy is congruent to a Chern--Simons
invariant of the boundary. $\Diamond$
\label{spektrum1}
\end{theorem}

\begin{remark}\rm It is clear that the above result depends only on 
the conformal class of the metric. One finds a similar 
energy identity for manifolds with conformally cylindrical 
ends \cite{mor-mro-rub} (including ALE spaces, in accordance with the 
energies of explicit instanton solutions of \cite{ete-hau2} and 
\cite{ete-hau3}) and for manifolds conformally of the form 
$\C\times\Sigma$ as in \cite{weh}. We expect that the validity of 
identities of this kind is more general. 

Taking into account the second part of Lemma \ref{cslemma} and the
fact that the character variety of a compact three-manifold has
finitely many connected components, we conclude that {\it the energy spectrum
of smooth, finite energy, self-dual connections over ALF spaces which 
satisfy the weak holonomy condition and decay rapidly is discrete}.

For irreducible instantons, imposing rapid decay is 
necessary for having discrete energies. For instance, in principle the 
energy formula \cite[Equation 2.32]{nye} provides a continuous 
energy spectrum for calorons and calorons of fractional energy are
known to exist (cf. e.g. \cite{der}). But slowly 
decaying reducible instantons still can have discrete spectrum; this is the 
case e.g. over the Schwarzschild space, cf. Section \ref{secfour}.  
 
Alternatively, instead of the rapid decay condition (\ref{lecsenges}), one
could also impose the possibly weaker but less natural condition that
the gauge invariant limit
\[\lim_{\varepsilon\rightarrow 0} \vert\tau_{\partial
M_\varepsilon}(A_\varepsilon ) -\tau_{\partial
M_\varepsilon}(\Gamma_\varepsilon )\vert = \mu\]
exists. Then the identity of Theorem \ref{spektrum1} would become:
\[\Vert F_A \Vert^2_{L^2(M,g)} \equiv
\tau_N (\Gamma_\infty ) + \mu ~~~~~ \mbox{{\rm mod} $\Z$}.\]
By analogy with the energy formula for calorons \cite[Equation 2.32]{nye},
we believe that the extra term $\mu$ is related to the {\it overall
magnetic charge} of an instanton while the modified energy formula would
represent the decomposition of the energy into ``electric''
(i.e., Chern--Simons) and ``magnetic'' (i.e., proportional to the
$\mu$-term) contributions.

In general, proving the existence of limits for the Chern--Simons 
functional assuming only the finiteness of the energy of the connection 
$\nabla_A$ is a very hard analytical problem, cf. \cite{mor-mro-rub,weh}. 
Therefore our rapid decay condition is a simple and natural condition 
which allows us to explicitly compute the limit of the Chern--Simons 
functional in our situation.

Another example illustrates that the weak holonomy condition is also
essential in Theorem \ref{spektrum1}. Consider $\R^4$, equipped with 
the Taub--NUT metric. This geometry admits a smooth $L^2$ harmonic 
2-form which can be identified with the curvature $F_B$ of a self-dual, 
rapidly decaying U$(1)$ connection $\nabla_B$ as in \cite{ete-hau2}; hence 
$\nabla_B\oplus\nabla^{-1}_B$ is a smooth, self-dual, rapidly 
decaying, reducible SU$(2)$ connection. We know that $H^2(\R^4 ,\Z )=0$ 
hence $\nabla_B$ lives on a trivial line bundle consequently it can be 
rescaled by an arbitrary constant like $B\mapsto cB$ without destroying 
its self-duality and finite energy. But the smooth, self-dual family
$\nabla_{cB}$ has continuous energy proportional to $c^2$. This strange
phenomenon also appears over the multi-Taub--NUT spaces, although they are
no more topologically trivial, cf. Section \ref{secfour} for more details.

From our holonomy viewpoint, this anomaly can be understood as follows.
Let $i: U^*_\varepsilon\subset W$ be an elementary neighbourhood as in
(\ref{elemikornyezet}) with the induced map $i_*:\pi_1(U^*_\varepsilon
)\rightarrow\pi_1(W)$. On the one hand we have $\pi_1(U^*_\varepsilon )\cong\Z$
as usual. On the other hand for the Taub--NUT space the asymptotical
topology is $W\cong S^3\times\R^+$ hence $\pi_1(W)\cong 1$ consequently
$i_*$ has a non-trivial kernel. However for a generic $c$ the connection
$\nabla_{cB}\vert_{U^*_\varepsilon}$ has non-trivial local holonomy 
$m\not= 0$ hence it does not obey the weak holonomy condition 
(\ref{holonomia}) therefore Theorem \ref{spektrum1} fails in this case.

The flat $\R^3\times S^1$ space has contrary behaviour to the
multi-Taub--NUT geometries. In this case we find $W\cong S^2\times
S^1\times\R^+$ for the end consequently $\pi_1(W)\cong\Z$ and the map $i_*$
is an obvious isomorphism. Hence the weak holonomy condition is always
obeyed. The character variety of the boundary is $\chi (S^2\times
S^1)\cong [0,1)$ hence connected. Referring to the second part of Lemma
\ref{cslemma} we conclude then that the energy of any smooth, self-dual,
rapidly decaying connection over the flat $\R^3\times S^1$ must be a 
non-negative integer in accordance with the known explicit solutions 
\cite{bru-nog-baa}. The case of the Schwarzschild space is simlar, cf. 
Section \ref{secfour}.
\end{remark}

\noindent These observations lead us to a more transparent form of Theorem
\ref{spektrum1} by replacing the weak holonomy condition with a simple,
sufficient topological criterion, which amounts to a straightforward
re-formulation of Definition \ref{whc}.

\begin{theorem} Let $(M,g)$ be an ALF space with an end
$W\cong N\times\R^+$ as before and, referring to the fibration 
(\ref{perem}), assume $N$ is an arbitrary circle bundle over 
$B_\infty\not\cong S^2,\:\R P^2$ or is a trivial circle bundle over $S^2$ 
or $\R P^2$. 

Then if $E$ is the (trivial) SU$(2)$ vector bundle over $M$ with a smooth, 
finite energy connection $\nabla_A$ then it satisfies the weak holonomy 
condition (\ref{holonomia}). 

Moreover if $\nabla_A$ is self-dual and decays rapidly 
as in (\ref{lecsenges}) then its energy is congruent to one of the 
Chern--Simons invariants of the boundary $N$.

In addition if the character variety $\chi (N)$ is connected, then the 
energy of any smooth, self-dual, rapidly decaying connection must be a 
non-negative integer.
\label{spektrum2}
\end{theorem}

\noindent{\it Proof.} Consider an elementary neighbourhood $i:
U^*_\varepsilon\subset W$ as in (\ref{elemikornyezet}) and the induced map
$i_*:\pi_1(U^*_\varepsilon )\rightarrow\pi_1(W)$. If ${\rm
Ker}\:i_*=\{0\}$ then $\nabla_A$ obeys (\ref{holonomia}). However one sees 
that $\pi_1(U^*_\varepsilon )\cong\pi_1(F)$ and $\pi_1(W)\cong\pi_1(N)$ 
hence the map $i_*$ fits well into the homotopy exact sequence 
\[\dots\longrightarrow\pi_2(B_\infty 
)\longrightarrow\pi_1(F)\stackrel{i_*}\longrightarrow
\pi_1(N)\longrightarrow\pi_1(B_\infty )\longrightarrow\dots\]
of the fibration (\ref{perem}). This segment shows that ${\rm 
Ker}\:i_*\not=\{ 0\}$ if and only if 
$N$ is either a non-trivial circle bundle over $S^2$ that is, $N\cong 
S^3/\Z_s$ a lens space of type $L(s,1)$, or a non-trivial circle bundle over 
$\R P^2$.

The last part is clear via the second part of Lemma \ref{cslemma}. $\Diamond$
\vspace{0.1in}

\noindent Finally we investigate the strong holonomy condition
(\ref{holonomia2}). As an important corollary of our construction we find

\begin{theorem}{\em (Sibner--Sibner, 1992; R\aa de, 1994)} Let $(M,g)$ be
an ALF space with an end $W$ as before and $E$ be the trivial SU$(2)$ 
vector bundle over $M$ with a smooth, finite energy connection
$\nabla_A$ and associated locally flat connection $\nabla_\Gamma$ on
$E\vert_W$. If and only if the strong holonomy condition
(\ref{holonomia2}) is satisfied then both $\nabla_A$ as
well as $\nabla_\Gamma$ as a flat connection, extend smoothly over the
whole $X$, the Hausel--Hunsicker--Mazzeo compactification of $(M,g)$. That
is, there exist bundles $\tilde{E}$ and $\tilde{E}_0\cong X\times\C^2$ 
over $X$ such that $\tilde{E}\vert_M\cong E$ and the connection $\nabla_A$
extends smoothly over $\tilde{E}$; in the same fashion
$\tilde{E}_0\vert_W\cong E\vert_W$ and $\nabla_\Gamma$ extends
smoothly as a flat connection over $\tilde{E}_0$.
\label{kiterjesztes}
\end{theorem}

\noindent{\it Proof.} The restriction of the embedding
$M\subset X$ gives $U^*_\varepsilon\subset U_\varepsilon$ and this
later space is contractible. Consequently if $i:
U^*_\varepsilon\subset X$ is the embedding then for the induced map
$i_*:\pi_1(U^*_\varepsilon )\rightarrow \pi_1(X)$ we always have ${\rm
Ker}\:i_*=\pi_1(U^*_\varepsilon )$ hence the connections $\nabla_A$ and
$\nabla_\Gamma$ extend smoothly over $X$ via Theorem \ref{lokalissibner} 
if and only if the strong holonomy condition (\ref{holonomia2}) holds. In 
particular the extension of $\nabla_\Gamma$ is a flat connection. $\Diamond$
\begin{remark}\rm
If a finite energy self-dual connection satisfies 
the strong holonomy condition then its energy is integer via Theorem 
\ref{kiterjesztes} regardless its curvature decay. Consequently these 
instantons again have discrete energy spectrum. We may then ask ourselves about 
the relationship between the strong holonomy condition on the one hand and 
the weak holonomy condition 
imposed together with the rapid decay condition on the other hand.
\end{remark}

%%%%%%%%%%%%%%%%%%%%%%%%%%%%%%%%%%%%%%%%%%%%%%%%%%%%%%%%
%%%%%%%%%%%%%%%%%%%%%%%%%%%%%%%%%%%%%%%%%%%%%%%%%%%%%%%%

\section{The moduli space}\label{secthree}

In this Section we are going to prove that the moduli spaces of framed SU(2)
instantons over ALF manifolds form smooth, finite dimensional manifolds,
whenever non-empty. The argument will go along the by now familiar
lines consisting of three steps: (i) Compute the dimension of the space of 
infinitesimal deformations of an irreducible, rapidly decaying self-dual 
connection, satisfying the weak holonomy condition, using a variant of the 
Gromov--Lawson relative index theorem \cite{gro-law} and a vanishing theorem; 
(ii) Use the Banach space inverse and implicit function theorems to integrate 
the infinitesimal deformations and obtain a local moduli space; (iii) Show 
that local moduli spaces give local coordinates on the global moduli space 
and that this global space is a Hausdorff manifold. We will carry out the 
calculations in detail for step (i) while just sketch (ii) and (iii) and refer 
the reader to the classical paper \cite{ati-hit-sin}.

Let $(M,g)$ be an ALF space with a single end $W$ as in Section \ref{sectwo}.
Consider a trivial SU$(2)$ bundle $E$ over $M$ with a smooth, 
irreducible, self-dual, finite energy connection $\nabla_A$ on it. 
By smoothness we mean that the connection 1-form is smooth in any smooth 
trivialization of $E$. In addition suppose $\nabla_A$ satisfies the weak 
holonomy condition (\ref{holonomia}) as well as decays rapidly in the 
sense of (\ref{lecsenges}). Then by Theorem \ref{spektrum1} its energy is 
determined by a Chern--Simons invariant. We will assume that this energy 
$e:=\Vert F_A\Vert^2_{L^2(M,g)}$ is {\it fixed}.
 
Consider the associated flat connection $\nabla_\Gamma$ with 
holonomy $m\in [0,1)$ as in Theorem \ref{spektrum1}. Extend 
$\nabla_\Gamma$ over the whole $E$ and 
continue to denote it by $\nabla_\Gamma$. Take the smooth gauge 
(\ref{globalismertek}) on the neck. Since $E$ is trivial, we can extend 
this gauge smoothly over the whole $M$ and can write $\nabla_A =\dd+A$ and 
$\nabla_\Gamma =\dd +\Gamma$ for some smooth connection 1-forms $A$ 
and $\Gamma$ well defined over the whole $M$. We also fix this gauge once 
and for all in our forthcoming calculations.
In particular the asymptotics of $\nabla_A$ is also {\it fixed} and is 
given by $\Gamma$. The connection $\nabla_\Gamma$, the usual Killing form 
on End$E$ and the rescaled metric $\tilde{g}$ are used to construct 
Sobolev spaces over various subsets of $X$ with respect to the norm 
(\ref{normak}). Both the energy $e$ and the asymptotics $\Gamma$ are 
preserved under gauge transformations which tend to the identity with 
vanishing first derivatives everywhere in infinity. We suppose ${\rm 
Aut}E\subset{\rm End}E$ and define the gauge group to be the completion 
\[{\cal G}_E:=\overline{\{\gamma -1\in C_0^\infty ({\rm 
End}E)\:\vert\:\Vert\gamma -1\Vert_{L^2_{j+2,\Gamma}(M)}<\infty ,\gamma\in 
\mbox{$C^\infty ({\rm Aut}E)$ a.e.}\}}\] 
and the gauge equivalence class of $\nabla_A$ under ${\cal G}_E$ is 
denoted by $[\nabla_A]$. Then we are seeking the virtual dimension of the 
framed moduli space ${\cal M}(e,\:\Gamma )$ of all such connections given up 
to these specified gauge transformations.

Writing $\dd_A: C^\infty (\Lambda^pM\otimes{\rm End}E)\rightarrow 
C^\infty (\Lambda^{p+1}M\otimes{\rm End}E)$ consider the usual deformation 
complex \[L^2_{j+2,\Gamma}(\Lambda^0M\otimes{\rm End 
}E)\stackrel{\dd_A}{\longrightarrow}L^2_{j+1,\Gamma}(\Lambda^1M\otimes{\rm
End}E) \stackrel{\dd^-_A}{\longrightarrow}L^2_{j,
\Gamma}(\Lambda^-M\otimes{\rm End}E)\]
where $\dd_A^-$ refers to the induced connection composed with the
projection onto the anti-self-dual side. Our first step is to check that
the Betti numbers $h^0,h^1,h^-$ of this complex, given by 
$h^0=\dim H^0({\rm End}E)$, etc., are finite. We therefore introduce an 
elliptic operator
\begin{equation}
\delta^*_A: L^2_{j+1,\Gamma}(\Lambda^1M\otimes{\rm 
End}E)\longrightarrow 
L^2_{j,\Gamma}((\Lambda^0M\oplus\Lambda^-M)\otimes{\rm End}E),
\label{komplexus}
\end{equation}
the so-called {\it deformation operator} 
$\delta^*_A:=\dd_A^*\oplus\dd_A^-$,
which is a conformally invariant first order elliptic operator
over $(M,\tilde{g})$ hence $(M,g)$. Here $\dd_A^*$ is the formal 
$L^2$ adjoint of $\dd_A$. We will demonstrate that $\delta^*_A$ is 
Fredholm, so it follows that $h^1=\dim{\rm Ker}\:\delta^*_A$ and 
$h^0+h^-=\dim{\rm Coker}\:\delta^*_A$ are finite. 

Pick up the trivial flat SU$(2)$ connection $\nabla_{\Theta}$
on $E$; it satisfies the strong holonomy condition
(\ref{holonomia2}), hence it extends smoothly over $X$ to an operator 
$\nabla_{\tilde{\Theta}}$ by Theorem \ref{kiterjesztes}.
Using the regularized metric $\tilde{g}_\varepsilon$ of (\ref{metrika})
it gives rise to an induced elliptic operator over the compact space 
$(X,\tilde{g}_\varepsilon )$ as
\[\delta^*_{\varepsilon, \tilde{\Theta}}: C_0^\infty (\Lambda^1X\otimes{\rm
End}\tilde{E}_0)\longrightarrow
C_0^\infty ((\Lambda^0X\oplus\Lambda^-X)\otimes{\rm End}\tilde{E}_0).\]
Consequently $\delta^*_{\varepsilon, \tilde{\Theta}}$ hence its 
restrictions $\delta^*_{\varepsilon ,\tilde{\Theta}}\vert_M$ and 
$\delta^*_{\varepsilon ,\tilde{\Theta}}\vert_W$ are Fredholm with respect to
any Sobolev completion. We construct a particular completion as follows. 
The smooth extended connection $\nabla_\Gamma$ on $E$ can be
extended further over $X$ to a connection $\nabla_{\tilde\Gamma}$ such
that it gives rise to an SU$(2)$ Sobolev connection on the
trivial bundle $\tilde{E}_0$. Consider a completion like
\begin{equation}
\delta^*_{\varepsilon, 
\tilde{\Theta}}:L^2_{j+1,\tilde{\Gamma}}(\Lambda^1X\otimes{\rm
End}\tilde{E}_0)\longrightarrow
L^2_{j,\tilde{\Gamma}}((\Lambda^0X\oplus\Lambda^-X)\otimes{\rm 
End}\tilde{E}_0).
\label{laposkomplexus}
\end{equation}
The operators in (\ref{komplexus}) and (\ref{laposkomplexus}) give rise to
restrictions. The self-dual connection yields
\[\delta^*_A\vert_W: L^2_{j+1,\tilde{\Gamma}}(\Lambda^1W\otimes{\rm
End}E\vert_W)\longrightarrow
L^2_{j,\tilde{\Gamma}}((\Lambda^0W\oplus\Lambda^-W)\otimes {\rm 
End}E\vert_W)\]
while the trivial connection gives
\[\delta^*_{\varepsilon ,\tilde{\Theta}}\vert_W:
L^2_{j+1,\tilde{\Gamma}}(\Lambda^1W\otimes{\rm
End}\tilde{E}_0\vert_W)\longrightarrow
L^2_{j,\tilde{\Gamma}}((\Lambda^0W\oplus\Lambda^-W)\otimes
{\rm End}\tilde{E}_0\vert_W).\]
Notice that these operators actually act on isomorphic Sobolev spaces 
consequently comparing them makes sense. Working in these Sobolev spaces we 
claim that

\begin{lemma} The deformation operator $\delta^*_A$ of 
(\ref{komplexus}) with $j=0,1$ satisfies the operator norm inequality
\begin{equation}
\Vert (\delta^*_A-\delta^*_{\varepsilon , 
\tilde{\Theta}})\vert_{V^*_\varepsilon}\Vert\leq 3\cdot 2^j\left( c_3\Vert 
F_A\Vert_{L^2(V_\varepsilon )}+c_4m\varepsilon^5\right)
\label{opinorm}
\end{equation}
where $m\in [0,1)$ is the holonomy of $\nabla_A$ and 
$c_4=c_4(\tilde{g}\vert_{V_\varepsilon},\gamma_\alpha ,\dd\gamma_\alpha 
)>0$ is a constant depending only on the metric and the gauge transformations 
used in (\ref{ujmerce}). Consequently $\delta^*_A$ is a Fredholom operator 
over $(M,g)$.
\label{fredholm}
\end{lemma}

\noindent{\it Proof.} Consider the restriction of the operators constructed 
above to the neck $W\cong V^*_\varepsilon$ and calculate the operator norm 
of their difference with an 
$a\in L^2_{j+1,\Gamma}(\Lambda^1M\otimes{\rm End}E)$ as follows:
\[\Vert (\delta^*_A-\delta^*_{\varepsilon 
,\tilde{\Theta}})\vert_{V^*_\varepsilon}\Vert 
=\sup\limits_{a\not =0}\frac{ \Vert (\delta^*_A-\delta^*_{\varepsilon 
,\tilde{\Theta}})a\Vert_{L^2_{j,\Gamma}(V^*_\varepsilon 
)}}{\Vert a\Vert_{L^2_{j+1,\Gamma}(V^*_\varepsilon )}}.\]
By assumption $\nabla_A$ satisfies
the weak holonomy condition (\ref{holonomia}) hence the connection 
$\nabla_\Gamma$ is flat on $E\vert_{V^*_\varepsilon}$ and it determines a 
deformation operator $\delta^*_{\tilde{\Gamma}}\vert_{V^*_\varepsilon}$ 
over $(V^*_\varepsilon ,\tilde{g}\vert_{V^*_\varepsilon})$. We use the triangle 
inequality: 
\[\Vert (\delta^*_A -\delta^*_{\varepsilon,
\tilde{\Theta}})\vert_{V^*_\varepsilon}\Vert\leq\Vert 
(\delta^*_A-\delta^*_{\tilde{\Gamma}} )\vert_
{V^*_\varepsilon}\Vert +
\Vert (\delta^*_{\tilde{\Gamma}} -\delta^*_{\varepsilon
,\tilde{\Theta}})\vert_{V^*_\varepsilon}\Vert .\]
Referring to the global gauge (\ref{globalismertek}) and the 
metric $\tilde{g}$ we have $\delta^*_A a 
=(\delta +\dd^-)a +A^*a +(A\wedge a +a\wedge A )^-$ 
and the same for $\delta^*_{\tilde{\Gamma}}$ consequently
\[(\delta^*_A -\delta^*_{\tilde{\Gamma}} )\vert_
{V^*_\varepsilon} =(A_{V^*_\varepsilon}-\Gamma_{V^*_\varepsilon} 
)^*+((A_{V^*_\varepsilon}-\Gamma_{V^*_\varepsilon})\wedge\cdot
+\cdot\wedge (A_{V^*_\varepsilon}-\Gamma_{V^*_\varepsilon}))^-.\]
Taking $j=0,1$ and combining this with (\ref{globalissibner}) we find an 
estimate for the first term like
\[\Vert (\delta^*_A-\delta^*_{\tilde{\Gamma }})\vert_{V^*_\varepsilon}
\Vert\leq 3\cdot 2^j\Vert A_{V^*_\varepsilon}-\Gamma_{V^*_\varepsilon} 
\Vert_{L^2_{j,\Gamma}(V^*_\varepsilon 
)}\leq 3\cdot 2^jc_3\Vert F_A\Vert_{L^2(V_\varepsilon )}.\]
Regarding the second term, the trivial flat connection
$\nabla_{\Theta}$ on $E\vert_{V^*_\varepsilon}$ satisfies the
strong holonomy condition (\ref{holonomia2}). In the
gauge (\ref{globalismertek}) we use, we may suppose
simply $\Theta\vert_{V^*_\varepsilon}=0$
consequently neither data from $\nabla_{\tilde{\Theta}}$ nor the perturbed
metric $\tilde{g}_\varepsilon$ influence this term. Now take a partition
of the end into elementary neighbourhoods (\ref{elemikornyezet})
and use the associated simple constant gauges (\ref{lokalismertek}) then
\[\Vert (\delta^*_{\tilde{\Gamma}} -\delta^*_{\varepsilon
,\tilde{\Theta}})\vert_{V^*_\varepsilon}\Vert
\leq 3\cdot 2^j\sum\limits_\alpha 
\Vert\gamma^{-1}_\alpha\Gamma_m\gamma_\alpha
\Vert_{L^2_{j,\Gamma}(U^*_{\varepsilon ,\alpha} 
)}\leq 3\cdot 2^jc_4m\varepsilon^5\] 
with some constant $c_4=c_4(\tilde{g}\vert_{V^*_\varepsilon},\gamma_\alpha 
,\dots ,\nabla^j_\Gamma\gamma_\alpha )>0$ via (\ref{terfogat}). Putting 
these together we get (\ref{opinorm}).

Taking into account that the right hand side of (\ref{opinorm}) is 
arbitrarily small we conclude that $\delta^*_A\vert_{V^*_\varepsilon}$ is 
Fredholm because so is $\delta^*_{\varepsilon 
,\tilde{\Theta}}\vert_{V^*_\varepsilon}$ and Fredholmness is an open property. 
Clearly, $\delta^*_A\vert_{M\setminus V^*_\varepsilon}$ is also Fredholm, 
because it is an elliptic operator
over a compact manifold. Therefore glueing the parametrices of these
operators together, one constructs a parametrix for $\delta^*_A$ over the
whole $M$ (see \cite{jar2} for an analogous construction). This shows that 
$\delta^*_A$ is Fredholm over the whole $(M,\tilde{g})$ hence $(M,g)$ as 
claimed. $\Diamond$
\vspace{0.1in}

\noindent We assert that $T_{[\nabla_A]}{\cal M}(e,\Gamma )\cong{\rm 
Ker}\:\delta^*_A$ consequently $h^1=\dim{\cal M}(e,\Gamma )$. Indeed, 
the anti-self-dual part of the curvature of a perturbed connection 
$\nabla_{A+a}$ is given by 
$F^-_{A+a}=\dd^-_Aa+(a\wedge a)^-$ with the perturbation $a\in
L^2_{j+1,\Gamma}(\Lambda^1M\otimes{\rm End}E)$. Therefore if $a\in{\rm 
Ker}\:\delta^*_A$ then both self-duality and the energy of $\nabla_{A+a}$ 
are preserved {\it infinitesimally} i.e., in first order. The perturbation
vanishes everywhere in infinity hence the asymptotics given by $\Gamma$ is 
also unchanged; in particular $\nabla_{A+a}$ continues to obey the weak 
holonomy condition. 

Incidentally, we note however that {\it locally} hence also {\it 
globally}, some care is needed when one perturbs a connection in our moduli 
spaces. For a perturbation with $a\in 
L^2_{j+1,\Gamma}(\Lambda^1M\otimes{\rm End}E)$ the asymptotics and in 
particular the weak holonomy condition are obeyed as we have seen. 
But concerning the energy, by repeating the 
calculation of Section \ref{sectwo} again, we obtain 
\[\frac{1}{\varepsilon}\int\limits_0^\varepsilon\vert\tau_{\partial
M_\rho}(A_\rho +a_\rho )-\tau_{\partial M_\rho}(\Gamma_\rho
)\vert^2\dd\rho\leq\]
\[\leq\frac{1}{\varepsilon}\int\limits_0^\varepsilon\left(\vert
\tau_{\partial M_\rho}(A_\rho +a_\rho )-\tau_{\partial
M_\rho}(A_\rho )\vert +\vert\tau_{\partial
M_\rho}(A_\rho )-\tau_{\partial M_\rho}(\Gamma_\rho
)\vert\right)^2\dd\rho\leq\]
\[\leq\left(\frac{c_1}{\sqrt{\varepsilon}}\Vert 
a\Vert_{L^2_{1,\Gamma}(V_\varepsilon
)}+\frac{c_1c_3}{\sqrt{\varepsilon}}\Vert
F_A\Vert_{L^2(V_\varepsilon )}\right)^2.\]
If the last line tends to zero as $\varepsilon\rightarrow 0$ then the 
perturbed connection has the same limit as in (\ref{magnes}).
Consequently, we find that the energy is also unchanged. Since the original 
connection decays rapidly in the sense of (\ref{lecsenges}) if the 
perturbations also decays rapidly i.e., if 
$\lim\limits_{\varepsilon\rightarrow 0}\varepsilon^{-\frac{1}{2}}\Vert 
a\Vert_{L^2_{1,\Gamma}(V_\varepsilon )}=0$ then the energy is preseved by 
local (i.e., small but finite) perturbations as well. 

It is therefore convenient to introduce {weighted Sobolev spaces} with 
weight $\delta =\frac{1}{2}$ and to say that $a$ and $\nabla_A$ {\it decay 
rapidly} if $a\in L^2_{\frac{1}{2},j+1,\Gamma}(\Lambda^1M\otimes{\rm End}E)$ 
and $F_A\in L^2_{\frac{1}{2},j,\Gamma}(\Lambda^2M\otimes{\rm End}E)$, 
respectively. These are gauge 
invariant conditions and for $\nabla_A$ with $j=0$ it is equivalent to 
(\ref{lecsenges}). In this framework the rough estimate
\[\Vert F_{A+a}\Vert_{L^2_{\frac{1}{2},j,\Gamma}(M)}=\Vert
F_A+\dd_Aa+a\wedge a\Vert_{L^2_{\frac{1}{2},j,\Gamma}(M
)}\leq\Vert F_A\Vert_{L^2_{\frac{1}{2},j,\Gamma}(M
)}+\Vert\dd_Aa\Vert_{L^2_{\frac{1}{2},j,\Gamma}(M
)}+c_5\Vert a\Vert^2_{L^2_{\frac{1}{2},j+1,\Gamma}(M)}\]
with some constant $c_5=c_5(\tilde{g})>0$ implies that $\nabla_{A+a}$ also 
decays rapidly. 

Let ${\cal A}_E$ denote the affine space of rapidly 
decaying SU$(2)$ connections on $E$ as well as ${\cal F}^-_E$ the vector 
space of the anti-self-dual parts of their curvatures. Then take a complex of 
punctured spaces, the global version of the deformation complex above:
\[({\cal G}_E, 1)\stackrel{f_A}{\longrightarrow}({\cal A}_E,
\nabla_A)\stackrel{\varrho^-_A}{\longrightarrow}({\cal F}^-_E, 0).\]
We have $({\cal A}_E,\nabla_A)\cong 
L^2_{\frac{1}{2},j+1,\Gamma}(\Lambda^1M\otimes{\rm End}E)$ and $({\cal 
F}^-_E, 0)\cong L^2_{\frac{1}{2},j,\Gamma}(\Lambda^-M\otimes{\rm End}E)$.
The global gauge fixing map is defined as $f_A(\gamma
):=\gamma^{-1}\nabla_A\gamma -\nabla_A$ while $\varrho^-_A(a):=F^-_{A+a}
=\dd^-_Aa+(a\wedge a)^-$. If $\nabla_A$ is irreducible then ${\rm 
Ker}\:f_A\cong 1$. One easily shows that if both $\gamma -1$ 
and $a$ are pointwise small then $f_A(\gamma )=\nabla_A(\log\gamma )$ and 
its formal $L^2$ adjoint satisfying
\[(f_A(\gamma ),a)_{L^2({\cal A}_E, 
\nabla_A)}=(\gamma ,f^*_A(a))_{L^2({\cal G}_E,1)}\]
looks like $f^*_A(a)=\exp (\nabla_A^*a)$ if 
\[(\gamma ,\beta )_{L^2({\cal G}_E,1)}
:=-\int\limits_M{\rm tr}(\log\gamma\log\beta )*_{\tilde{g}}1,\]
defined in a neighbourhood of $1\in{\cal G}_E$. If $\Delta$ denotes the 
diagonal of $({\cal A}_E,\nabla_A)\times ({\cal A}_E,\nabla_A)$ we obtain a model for 
$O_{[\nabla_A]}\subset {\cal M}(e,\Gamma )$, the vicinity of $[\nabla_A]$, 
as follows:
\[O_{[\nabla_A]}\cong{\rm Ker}(f^*_A\times \varrho^-_A)\vert_{\Delta}\subset 
({\cal A}_E,\nabla_A).\]
The derivative of $(f^*_A\times \varrho^-_A)\vert_{\Delta}$ at $a$ is 
$\delta^*_A+( a\wedge\cdot +\cdot\wedge a)^-$ which is a Fredholm operator hence 
$O_{[\nabla_A]}$ is smooth and finite dimensional. These local models 
match together and prove that the moduli space is indeed a smooth manifold 
of dimension $h^1$.

We return to the calculation of $h^1$. We will calculate the index of 
$\delta^*_A$ in (\ref{komplexus}) by referring to a relative index 
theorem. This provides us with the alternating sum $-h^0+h^1-h^-$ and 
then we show that $h^0=h^-=0$ via a vanishing theorem.

First we proceed to the calculation of the index of $\delta^*_A$. This 
will be carried out by a variant of the Gromov--Lawson relative index theorem 
\cite{gro-law}, which we will now explain. First, let us introduce some 
notation. For any elliptic Fredholm operator
$P$, let ${\rm Index}_aP$ denote its analytical index, i.e., ${\rm
Index}_aP=\dim{\rm Ker}P-\dim{\rm Coker}P$. If this $P$ is defined over a
compact manifold, ${\rm Index}_aP$ is given by a topological formula as in
the Atiyah--Singer index theorem, which we denote by ${\rm Index}_tP$,
the topological index of $P$. The following theorem will be proved in the 
Appendix (cf. \cite{jar1}):

\begin{theorem}
Let $(M,g)$ be a complete Riemannian manifold, and let $X$ be some smooth
compactification of $M$. Let also $D_1: L^2_{j+1,\Gamma_1}(F_1)
\rightarrow L^2_{j,\Gamma '_1}(F_1')$ and $D_0: L^2_{j+1,\Gamma_0}(F_0)
\rightarrow L^2_{j,\Gamma '_0}(F_0')$ be 
two first order, elliptic Fredholm operators defined on complex vector 
bundles $F_1,F_1'$ and $F_0,F_0'$ over $M$ with fixed Sobolev connections
$\nabla_{\Gamma_1},\nabla_{\Gamma '_1}$ and 
$\nabla_{\Gamma_0},\nabla_{\Gamma '_0}$, respectively.
Assume that given $\kappa>0$, there is a compact subset
$K\subset M$ such that, for $W=M\setminus K$, the following hold:
\begin{itemize}

\item[(i)] There are bundle isomorphisms $\phi: F_1\vert_W\cong
F_0\vert_W$ and $\phi ': F_1'\vert_W\cong F_0'\vert_W$;

\item[(ii)] The operators asymptotically agree, that is, in some
operator norm $\Vert (D_1-D_0)\vert_W\Vert <\kappa$.
\end{itemize}
If arbitrary elliptic extensions $\tilde{D}_1$ and 
$\tilde{D}_0$ of $D_1$ and $D_0$ to $X$ exist, then we have
\[{\rm Index}_aD_1-{\rm Index}_aD_0={\rm Index}_t\tilde{D}_1-{\rm
Index}_t\tilde{D}_0\]
for the difference of the analytical indices.
\label{index}
\end{theorem}

\noindent In the case at hand, $F_1=F_0={\Lambda^1M\otimes\rm 
End}E\otimes\C$ and
$F_1'=F_0'=(\Lambda^0M\oplus\Lambda^-M)\otimes{\rm End}E\otimes\C$;
then the complexification of $\delta^*_A$, which by Lemma 
\ref{fredholm} is a Fredholm operator, plays the role of $D_1$ and via 
(\ref{opinorm}) it asymptotically agrees with the complexified 
$\delta^*_{\varepsilon, \tilde{\Theta}}\vert_M$, also Fredholm by 
construction, which replaces $D_0$.

To find the operators which correspond to $\tilde{D}_1$ and $\tilde{D}_0$,
we proceed as follows. Remember that in general $\nabla_A$ does not extend
over $X$ (cf. Theorem \ref{kiterjesztes}). Let $\tilde{E}$ be
the unique vector bundle over $X$ constructed as follows. 
Since $\tilde{E}\vert_{X\setminus V_{2\varepsilon}}\cong 
E\vert_{M\setminus 
V^*_{2\varepsilon}}$ and $\tilde{E}\vert_{V_\varepsilon}\cong 
V_\varepsilon\times\C^2$, this bundle is uniquely determined by the glued 
connection $\nabla_{\tilde{A}}=(1-f_\varepsilon 
)\nabla_A+f_\varepsilon\nabla_{\tilde{\Theta}}$ on $\tilde{E}$ where 
$f_\varepsilon$ is taken from (\ref{metrika}). We can
construct the associated operator $\delta^*_{\varepsilon, \tilde{A}}$ over
${\rm End}\tilde{E}$ with respect to the metric (\ref{metrika}) whose
complexification will correspond to $\tilde{D}_1$; the operator 
$\tilde{D}_0$ is given by the complexification of the operator
$\delta^*_{\varepsilon ,\tilde{\Theta}}$ on the trivial bundle ${\rm 
End}\tilde{E}_0$, which we have already constructed.

The right hand side of the relative index formula in Theorem \ref{index}
is given by
\[{\rm Index}_t(\delta^*_{\varepsilon ,\tilde{A}})-{\rm
Index}_t(\delta^*_{\varepsilon
,\tilde{\Theta}})=8k-3(1-b^1(X)+b^-(X))+3(1-b^1(X)+b^-(X))=8k\]
where $k=\Vert F_{\tilde{A}}\Vert^2_{L^2(X)}$ is the second
Chern number of the extended bundle $\tilde{E}$.
Notice that this number might be different from the energy $e=\Vert
F_A\Vert^2_{L^2(M,g)}$ of the original connection. We only know a priori
that $k\leq e$. However we claim that
\begin{lemma}
Using the notation of Theorem \ref{spektrum1} and, in the same fashion, if 
$\nabla_\Theta =\dd +\Theta_W$ on $E\vert_W$ then letting $\Theta_\infty
:=\lim\limits_{r\rightarrow\infty}\Theta_W\vert_{N\times\{r\}}$, in any 
smooth gauge 
\begin{equation}
k=e+\tau_N(\Theta_\infty )-\tau_N(\Gamma_\infty )
\label{e-formula}
\end{equation}
holds. Notice that this expression is gauge invariant.
\label{hiba}
\end{lemma}
\begin{remark}\rm
Of course in any practical application it is worth taking the gauge in 
which simply $\tau_N(\Theta_\infty )=0$.
\end{remark}

\noindent {\it Proof.} Using the gauge (\ref{globalismertek}) 
for instance and applying the Chern--Simons theorem for the restricted
energies to $M_\varepsilon$, we find
\[(k-e)\vert_{M_\varepsilon} =\tau_{\partial 
M_\varepsilon}(\tilde{A}_\varepsilon )-\tau_{\partial
M_\varepsilon}(A_\varepsilon ) +2\Vert 
F^-_{\tilde{A}}\Vert^2_{L^2(M_\varepsilon )}.\]
Therefore, since $\Vert F^-_{\tilde{A}}\Vert_{L^2(M_\varepsilon 
)}=\Vert F^-_{\tilde{A}}\Vert_{L^2(V_{2\varepsilon} )}$ and the 
Chern--Simons invariants converge as in (\ref{magnes}) by the rapid decay 
assumption, taking the limit one obtains
\[k-e=\tau_N(\Theta_0)-\tau_N(\Gamma_0)+2\lim\limits_{\varepsilon\rightarrow 
0}\Vert F^-_{\tilde{A}}\Vert^2_{L^2(V_{2\varepsilon})}
=\tau_N(\Theta_\infty )-\tau_N(\Gamma_\infty)+
2\lim\limits_{R\rightarrow\infty}\Vert F^-_{\tilde{A}}\Vert^2_
{L^2\left( V_{2R},g\vert_{V_{2R}}\right)}.\]
Consequently we have to demonstrate that 
$\lim\limits_{\varepsilon\rightarrow
0}\Vert F^-_{\tilde{A}}\Vert_{L^2(V_{2\varepsilon})}=0$. There is a 
decomposition of the glued curvature like $F_{\tilde{A}}=\Phi + \varphi$ 
with $\Phi :=(1-f_\varepsilon )F_A +f_\varepsilon F_\Theta$
and a perturbation term as follows: 
\[\varphi :=-\dd f_\varepsilon\wedge
(A_{V^*_{2\varepsilon}}-\Theta_{V^*_{2\varepsilon}}) -f_\varepsilon 
(1-f_\varepsilon )(A_{V^*_{2\varepsilon}}-\Theta_{V^*_{2\varepsilon}})\wedge 
(A_{V^*_{2\varepsilon}}-\Theta_{V^*_{2\varepsilon}}).\]
This shows that $F^-_{\tilde{A}}=\varphi^-$ consequently it is compactly 
supported in $V_{2\varepsilon}\setminus V_\varepsilon$. Moreover there is 
a constant $c_6=c_6(\dd f_\varepsilon ,\tilde{g}\vert_{V_{2\varepsilon}})>0$, 
independent of $\varepsilon$, such that $\Vert\dd 
f_\varepsilon\Vert_{L^2_{1,\Gamma}(V_{2\varepsilon})}\leq c_6$; as well as
$\vert f_\varepsilon (1-f_\varepsilon )\vert\leq\frac{1}{4}$ therefore, 
recalling the pattern of the proof of Lemma \ref{fredholm}, we tame 
$\varphi^-$ like 
\[\Vert\varphi^-\Vert_{L^2(V_{2\varepsilon})}\leq
\Vert\varphi\Vert_{L^2(V_{2\varepsilon})}\leq 
c_5c_6\Vert A_{V^*_{2\varepsilon}}-\Theta_{V^*_{2\varepsilon}}
\Vert_{L^2_{1,\Gamma}(V^*_{2\varepsilon})} 
+\frac{c_5}{4}\Vert A_{V^*_{2\varepsilon}}-\Theta_{V^*_{2\varepsilon}}
\Vert^2_{L^2_{1,\Gamma}(V^*_{2\varepsilon})}\leq\]
\[\leq c_5c_6\left( c_3\Vert 
F_A\Vert_{L^2(V_{2\varepsilon})}+c_4m(2\varepsilon )^5\right) +
\frac{c_5}{4}\left( c_3\Vert
F_A\Vert_{L^2(V_{2\varepsilon})}+c_4m(2\varepsilon )^5\right)^2.\] 
However we know that this last line can be kept as small as one likes 
providing the result. $\Diamond$
\vspace{0.1in}

\noindent Regarding the left hand side of Theorem \ref{index}, on the one 
hand we already know that
\[{\rm Index}_a\delta^*_A =-h^0+h^1-h^-=h^1=\dim{\cal M}(e,\Gamma )\]
by the promised vanishing theorem. On the other hand, since
${\rm End}E\otimes\C\cong M\times\C^3$, we find
\[{\rm Index}_a(\delta^*_{\varepsilon
,\tilde{\Theta}}\vert_M)=-3\left( b^0_{L^2}(M,\tilde{g}_\varepsilon\vert_M
)-b^1_{L^2}(M,\tilde{g}_\varepsilon\vert_M
)+b^-_{L^2}(M,\tilde{g}_\varepsilon\vert_M )\right)\]
where $b^i_{L^2}(M,\tilde{g}_\varepsilon\vert_M )$ is the $i$th 
$L^2$ Betti number and $b^-_{L^2}(M,\tilde{g}_\varepsilon\vert_M )$ is the 
dimension of the space of anti-self-dual finite energy 2-forms on the 
rescaled-regularized manifold $(M,\tilde{g}_\varepsilon\vert_M)$ 
i.e., this index is the truncated $L^2$ Euler characteristic of
$(M,\tilde{g}_\varepsilon\vert_M)$. We wish to cast this subtle invariant
into a more explicit form at the expense of imposing a further but
natural assumption on the spaces we work with.

\begin{lemma} Let $(M,g)$ be a Ricci flat ALF space, and let $X$ be its
compactification with induced orientation. Then one
has ${\rm Index}_a(\delta^*_{\varepsilon 
,\tilde{\Theta}}\vert_M)=-3b^-(X)$
where $b^-(X)$ denotes the rank of the negative definite part of the 
topological intersection form of $X$.
\label{L2lemma}
\end{lemma}

\noindent{\it Proof.} Exploiting the stability of the
index against small perturbations as well as the conformal invariance of 
the operator $\delta^*_{\varepsilon ,\tilde{\Theta}}\vert_M$, without 
changing the index we can replace the metric $\tilde{g}_\varepsilon$ with 
the original ALF metric $g$. Consequently we can write 
\[{\rm Index}_a(\delta^*_{\varepsilon ,\tilde{\Theta}}\vert_M)=
-3\left( b^0_{L^2}(M,g)-b^1_{L^2}(M,g)+b^-_{L^2}(M,g)\right)\]
for the index we are seeking.

Remember that this metric is complete and asymptotically looks like 
(\ref{aszimptotika}). This implies that $(M,g)$ has infinite volume hence 
a theorem of Yau \cite{yau} yields that $b^0_{L^2}(M,g)=0$. Moreover if we 
assume the curvature of $(M,g)$ not only satisfies (\ref{gorbulet}) but is
furthermore Ricci flat then a result of Dodziuk \cite{dod} shows
that in addition $b^1_{L^2}(M,g)=0$. Concerning $b^-_{L^2}(M,g)$ we use
the result of \cite{ete} (based on \cite[Corollary 7]{hau-hun-maz}) to 
observe that any finite energy anti-self-dual 2-form over $(M,g)$ extends 
smoothly as a (formally) anti-self-dual 2-form over $(X, \tilde{g})$ 
showing that $b^-_{L^2}(M,g)=b^-(X)$ as desired. $\Diamond$
\vspace{0.1in}

\noindent Finally we prove the vanishing of the numbers $h^0$ and $h^-$.
The proof is a combination of the standard method \cite{ati-hit-sin} and
a Witten-type vanishing result (\cite[Lemma 4.3]{par-tau}).
In the adjoint of the elliptic complex (\ref{komplexus}) we find
${\rm Index}_a\delta_A=\dim{\rm Ker}\:\delta_A-\dim{\rm
Coker}\:\delta_A=h^0-h^1+h^-$ hence proving ${\rm
Ker}\:\delta_A=\{0\}$ is equivalent to $h^0=h^-=0$.

\begin{lemma} Assume $(M,g)$ is an ALF space as defined in Section 2. 
If $\nabla_A$ is irreducible and
$\psi\in L^2_{j,\Gamma}((\Lambda^0M\oplus\Lambda^-M)\otimes{\rm 
End}E)$ satisfies $\delta_A\psi =0$ then $\psi =0$.
\label{eltunes}
\end{lemma}

\noindent{\it Proof.} Taking into account that ${\rm
Ker}\:\delta_A={\rm Ker}(\delta^*_A\delta_A)$ consists of smooth 
functions by elliptic regularity and is conformally invariant, we can use 
the usual Weitzenb\"ock formula with respect to the original ALF metric $g$ as 
follows:
\[\delta^*_A\delta_A=\dd_A^*\dd_A+W_g^-+\frac{s_g}{3}:\:\:\:
C_0^\infty ((\Lambda^0M\oplus\Lambda^-M)\otimes{\rm
End}E)\longrightarrow C_0^\infty
((\Lambda^0M\oplus\Lambda^-M)\otimes{\rm End}E).\]
In this formula $\dd_A: C_0^\infty
((\Lambda^0M\oplus\Lambda^-M)\otimes{\rm End}E)\rightarrow C_0^\infty
((\Lambda^1M\oplus\Lambda^3M)\otimes{\rm End}E)$
is the induced connection while $W_g^-+\frac{s_g}{3}$ acts only on the 
$\Lambda^-M$ summand as a symmetric, linear, algebraic map. This implies 
that if $\delta_A\psi =0$ such that
$\psi\in L^2_{j,\Gamma}(\Lambda^0M\otimes{\rm End}E)$ only 
or $W_g^-+\frac{s_g}{3}=0$ then $\nabla_A\psi =0$. If $\nabla_A$ is 
irreducible then both $\Lambda^0M\otimes{\rm End}E\otimes\C\cong 
S^0\Sigma^-\otimes{\rm End}E$ and $\Lambda^-M\otimes{\rm 
End}E\otimes\C\cong S^2\Sigma^-\otimes{\rm End}E$ are irreducible 
${\rm SU}(2)^-\times{\rm SU}(2)$ bundles hence $\psi =0$ follows.

Concerning the generic case, then as
before, we find $h^0=0$ by irreducibility consequently we can assume
$\psi\in L^2_{j, \Gamma}(\Lambda^-M\otimes{\rm End} E)$ only 
and $W_g^-+\frac{s_g}{3}\not=0$. Since ${\rm Ker}\:\delta_A$ consists of 
smooth functions it follows that if $\delta_A\psi =0$ then $\psi$ 
vanishes everywhere at infinity.

Let $\langle\:\cdot\:,\:\cdot\:\rangle$ be a pointwise SU$(2)$-invariant
scalar product on $\Lambda^-M\otimes{\rm EndE}$ and set
$\vert\psi\vert :=\langle\psi ,\psi\rangle^{\frac{1}{2}}$. Assume
$\delta_A\psi =0$ but $\psi\not =0$. Then
$\langle\dd^*_A\dd_A\psi ,\psi\rangle
=-\langle (W_g^-+\frac{s_g}{3})\psi ,\psi\rangle$ by the Weitzenb\"ock 
formula above. Combining this with the pointwise expression 
$\langle\dd_A^*\dd_A\psi ,\psi\rangle
=\vert\dd_A\psi\vert^2+\frac{1}{2}\triangle\vert\psi\vert^2$ and
applying $\frac{1}{2}\triangle\vert\psi\vert^2=
\vert\psi\vert\triangle\vert\psi\vert +\vert\dd\vert\psi\vert\vert^2$
as well as Kato's inequality
$\vert\dd\vert\psi\vert\vert\leq\vert\dd_A\psi\vert$,
valid away from the zero set of $\psi$, we obtain
\[2\:\frac{\vert\dd\vert\psi\vert\vert^2}{\vert\psi\vert^2}\leq\left\vert
W_g^-+\frac{s_g}{3}\right\vert 
-\frac{\triangle\vert\psi\vert}{\vert\psi\vert} 
=\left\vert W_g^-+\frac{s_g}{3}\right\vert
-\triangle\log\vert\psi\vert -\frac{\vert\dd\vert\psi\vert\vert^2}
{\vert\psi\vert^2}\]
that is,
\[3\:\vert\dd\log\vert\psi\vert\vert^2+\triangle\log\vert\psi\vert\leq
\left\vert W_g^-+\frac{s_g}{3}\right\vert .\]
Let $\lambda :\R^+\rightarrow W\cong N\times\R^+$ be a naturally
parameterized ray running toward infinity and let $f(r):=\log\vert\psi 
(\lambda (r))\vert$. Observe that $f$ is negative in the vicinity of a 
zero of $\psi$ or for large $r$'s. The last inequality then asymptotically 
cuts down along $\lambda$ to
\[3(f'(r))^2-\frac{2}{r}f'(r)-f''(r)\leq\frac{c_7}{r^3},\]
using the expansion of the Laplacian for a metric like
(\ref{aszimptotika}) and referring to the curvature decay
(\ref{gorbulet}) providing a constant $c_7=c_7(g)\geq 0$. This 
inequality yields $-c_7r^{-2}\leq (rf(r))''$. Integrating it we find
$c_7r^{-1}+a\leq (rf(r))'\leq rf'(r)$ showing $c_7r^{-2}+ar^{-1}\leq
f'(r)$ with some constant $a$ hence there is a constant $c_8\geq 0$ such 
that $c_8:=-\vert\inf\limits_{r\in\R^+}f'(r)\vert$. Integrating 
$c_7r^{-1}+a\leq (rf(r))'$ again we also obtain 
$c_7(\log r)r^{-1}+a+br^{-1}\leq f(r)\leq 0$, for some real constants 
$a,b$. 

Let $x_0\in M$ be such that $\psi (x_0)\not=0$; then by smoothness there
is another point $x\in M$ with this property such that $\vert x_0\vert <\vert
x\vert$. Integrating again the inequality
$-c_7r^{-2}\leq (rf(r))''$ twice from $x_0$ to $x$ along the ray
$\lambda$ connecting them we finally get
\[\vert\psi (x_0)\vert\leq\vert\psi
(x)\vert\exp\left( (c_7\vert x_0\vert^{-1}+c_8\vert x_0\vert )
(1-\vert x_0\vert\vert x\vert^{-1} )\right) .\]
Therefore either letting $x$ to be a zero of $\psi$ along $\lambda$ or, if 
no such point exists, taking the limit $\vert x\vert\rightarrow\infty$ we 
find that $\psi (x_0) =0$ as desired. $\Diamond$
\vspace{0.1in}

\noindent Finally, putting all of our findings together, we have arrived
at the following theorem:

\begin{theorem} Let $(M,g)$ be an ALF space with an end
$W\cong N\times\R^+$ as before. Assume furthermore that the metric is
Ricci flat. Consider a rank 2 complex SU$(2)$ vector bundle $E$ 
over $M$, necessarily trivial, and denote by ${\cal M}(e,\Gamma )$ the 
framed moduli space of smooth, irreducible, self-dual SU$(2)$ connections on 
$E$ satisfying the weak holonomy condition (\ref{holonomia}) and decaying 
rapidly in the sense of (\ref{lecsenges}) such that their energy $e<\infty$ is 
fixed and are asymptotic to a fixed smooth flat connection $\nabla_\Gamma$ on 
$E\vert_W$.

Then ${\cal M}(e,\Gamma )$ is either empty or a manifold of dimension
\[\dim{\cal M}(e,\Gamma )=8\left( e+\tau_N(\Theta_\infty
)-\tau_N(\Gamma_\infty )\right) -3b^-(X)\]
where $\nabla_{\Theta}$ is the trivial flat connection on $E\vert_W$ and 
$\tau_N$ is the Chern--Simons functional of the boundary 
while $X$ is the Hausel--Hunsicker--Mazzeo compactification of $M$ with 
induced orientation. $\Diamond$
\label{moduluster}
\end{theorem}
\begin{remark}\rm Of course, we get a dimension formula for anti-instantons by
replacing $b^-(X)$ with $b^+(X)$. Notice that our moduli spaces contain
framings, since we have a fixed flat connection and a gauge at infinity. 
The virtual dimension of the moduli space of unframed instantons is given 
by $\dim{\cal M}(e,\Gamma )-3$ which is is the number of effective free 
parameters.

A dimension formula in the presence of a magnetic term $\mu$ mentioned
in Section \ref{sectwo} is also easy to work out because in this case
(\ref{e-formula}) is simply replaced with $k=e+\tau_N(\Theta_\infty
)-\tau_N(\Gamma_\infty )-\mu$ and then this should be inserted into the
dimension formula of Theorem \ref{moduluster}.

Note also that our moduli spaces are naturally endowed with weighted $L^2$ 
metrics. An interesting problem is to investigate the properties of these 
metrics.
\end{remark}

%%%%%%%%%%%%%%%%%%%%%%%%%%%%%%%%%%%%%%%%%%%%%%%%%%%%%%%%%%%%%%%%%%%%%%%%
%%%%%%%%%%%%%%%%%%%%%%%%%%%%%%%%%%%%%%%%%%%%%%%%%%%%%%%%%%%%%%%%%%%%%%%%

\section{Case studies}\label{secfour}

In this Section we present some applications of Theorem
\ref{moduluster}. We will consider rapidly decaying instantons in the 
sense of (\ref{lecsenges}) over the 
flat $\R^3\times S^1$, the  multi-Taub--NUT geometries and the Riemannian 
Schwarzschild space. We also have the aim to enumerate the known 
Yang--Mills instantons over non-trivial ALF geometries. However we 
acknowledge that our list is surely incomplete, cf. e.g. \cite{ali-sac}.
\vspace{0.3in}

\centerline{\sc The flat space $\R^3\times S^1$} 
\vspace{0.2in}

\noindent This is the simplest ALF space hence instanton (or also called 
caloron i.e., instanton at finite temperature) theory over this space is 
well-known (cf. \cite{bru-baa,bru-nog-baa}). We claim that

\begin{theorem} Take $M=\R^3\times S^1$ with a fixed orientation and
put the natural flat metric onto it. Let $\nabla_A$ be a smooth, 
self-dual, rapidly decaying SU$(2)$ connection on a fixed rank 2
complex vector bundle $E$ over $M$. Then $\nabla_A$ satisfies the weak 
holonomy condition and has non-negative integer energy. Let ${\cal 
M}(e,\Gamma )$ denote the framed moduli space of these connections 
which are moreover irreducible as in Theorem \ref{moduluster}. Then
\[\dim {\cal M}(e,\Gamma ) =8e\]
and ${\cal M}(e,\Gamma )$ is not empty for all $e\in\N$.
\label{kaloron}
\end{theorem}

\noindent{\it Proof.} Since the metric is flat, the conditions of Theorem
\ref{moduluster} are satisfied. Furthermore, in the case at hand the
asymptotical topology of the space is $W\cong S^2\times S^1\times\R^+$ 
consequently $N\cong S^2\times S^1$ hence its 
character variety $\chi (S^2\times S^1)\cong [0,1)$ is connected. 
Theorem \ref{spektrum2} therefore guarantees 
that any smooth, self-dual, rapidly decaying connection has non-negative 
integer energy. In the gauge in which $\tau_{S^2\times 
S^1}(\Theta_\infty )=0$ we also get 
$\tau_{S^2\times S^1}(\Gamma_\infty )=0$. Moreover we find that 
$X\cong S^4$ for the Hausel--Hunsicker--Mazzeo 
compactification \cite{ete} yielding $b^-(X)=b^+(X)=0$; putting these data 
into the dimension formula in Theorem \ref{moduluster} we get the dimension 
as stated, in agreement with \cite{bru-baa}.

The moduli spaces ${\cal M}(e,\Gamma )$ are not empty for all $e\in\N$;
explicit solutions with arbitrary energy were constructed via a modified
ADHM construction in \cite{bru-baa,bru-nog-baa}. $\Diamond$

\vspace{0.3in}

\centerline{\sc The multi-Taub--NUT (or $A_k$ ALF, or ALF Gibbons--Hawking) 
spaces}
\vspace{0.2in}

\noindent The underlying manifold $M_V$ topologically can be understood as 
follows. There is a circle action on $M_V$ with $s$ distinct fixed points
$p_1,\dots,p_s\in M_V$, called NUTs. The quotient is $\R^3$ and we denote
the images of the fixed points  also by $p_1,\dots ,p_s\in \R^3$.
Then $U_V:=M_V\setminus \{ p_1,\dots, p_s\}$ is fibered over
$Z_V:=\R^3\setminus\{ p_1,\dots,p_s\}$ with $S^1$ fibers. The degree of
this circle bundle around each point $p_i$ is one.

The metric $g_V$ on $U_V$ looks like (cf. e.g. \cite[p. 363]{egu-gil-han})
\[\dd s^2=V(\dd x^2+\dd y^2+\dd z^2)+\frac{1}{V}(\dd\tau +\alpha )^2,\]
where $\tau\in (0,8\pi m]$ parameterizes the circles and
$x=(x,y,z)\in\R^3$; the smooth function $V: Z_V\rightarrow\R$ and the
1-form $\alpha\in C^\infty (\Lambda^1Z_V)$ are defined as follows:
\[V(x ,\tau )=V(x )=1+\sum\limits_{i=1}^s\frac{2m}{\vert x-p_i\vert}
,\:\:\:\:\:\dd\alpha =*_3\dd V.\]
Here $m>0$ is a fixed constant and $*_3$ refers to the Hodge-operation
with respect to the flat metric
on $\R^3$. We can see that the metric is independent of $\tau$ hence we
have a Killing field on $(M_V,g_V)$. This Killing field provides the above
mentioned U$(1)$-action. Furthermore it is possible
to show that, despite the apparent singularities in the NUTs, these
metrics extend analytically over the whole $M_V$ providing an ALF,
hyper--K\"ahler manifold. We also notice that the Killing field makes it
possible to write a particular K\"ahler-form in the hyper-K\"ahler family
as $\omega =\dd\beta$ where $\beta$ is a $1$-form of linear growth.

Then we can assert that

\begin{theorem} Let $(M_V, g_V)$ be a multi-Taub--NUT space with $s$ NUTs
and orientation induced by any of the complex structures in the
hyper--K\"ahler family. Consider the framed moduli space ${\cal
M}(e,\Gamma )$ of smooth, irreducible, rapidly decaying anti-self-dual 
connections satisfying the weak holonomy condition on a
fixed rank 2 complex SU$(2)$ vector bundle $E$ as in Theorem
\ref{moduluster}. Then ${\cal M}(e,\Gamma )$ is either empty or a manifold 
of dimension
\[\dim {\cal M}(e,\Gamma )=8\left( e+\tau_{L(s,-1)}(\Theta_\infty
)-\tau_{L(s,-1)}(\Gamma_\infty )\right)\] where $L(s,-1)$ is the 
oriented lens space representing the boundary of $M_V$. The moduli 
spaces are surely not empty for $\tau_{L(s,-1)}(\Theta_\infty
)=\tau_{L(s,-1)}(\Gamma_\infty )=0$ and $e=1,\dots, s$.
\label{taub}
\end{theorem}

\noindent{\it Proof.} This space is non-flat nevertheless its curvature 
satisfies the cubic curvature decay (\ref{gorbulet}) hence it is an 
ALF space in our sense. Since it is moreover hyper--K\"ahler, the 
conditions of Theorem \ref{moduluster} are satisfied. However this time the
asymptotic topology is $W\cong L(s,-1)\times\R^+$ therefore $N\cong 
L(s,-1)$ is a non-trivial circle bundle over $S^2$; consequently the 
weak holonomy condition (\ref{holonomia}) must be imposed. If the 
connection in 
addition decays rapidly as in (\ref{lecsenges}) then its energy is determined 
by a Chern--Simons invariant via Theorem \ref{spektrum1}. The character variety 
of the boundary lens space, $\chi (L(s,-1))$ is also non-connected if 
$s>1$ and each connected component has a non-trivial fractional Chern--Simons 
invariant which is calculable (cf., e.g. \cite{auc, kir-kla}). By the 
result in \cite{ete} the compactified space $X$ with its induced orientation is 
isomorphic to the connected sum of 
$s$ copies of $\overline{\C P}^2$'s therefore $b^+(X)=0$ and $b^-(X)=s$.
Inserting these into the dimension formula of Theorem \ref{moduluster} for
anti-self-dual connections we get the dimension.

Concerning non-emptiness, since lacking a general ADHM-like construction,
we may use a conformal rescaling method \cite{ete-hau2,ete-hau3}.
Take the natural orthonormal frame
\[\xi^0=\frac{1}{\sqrt{V}}(\dd\tau+\alpha),\:\:\:\:\:\xi^1=\sqrt{V}
\dd x,\:\:\:\:\:\xi^2=\sqrt{V}\dd y,\:\:\:\:\:\xi^3=\sqrt{V}\dd z\]
over $U_V$ and introduce the quaternion-valued 1-form $\boldsymbol{\xi}
:=\xi^0+\xi^1{\bf i}+\xi^2{\bf j}+\xi^3{\bf k}$. Moreover pick up the
non-negative function $f: U_V\rightarrow\R^+$ defined as
\[f(x):=\lambda_0+\sum\limits_{i=1}^s\frac{\lambda_i}{\vert x -p_i\vert}\]
with $\lambda_0,\lambda_1,\dots,\lambda_s$ being real
non-negative constants and also take the quaternion-valued $0$-form
\[{\bf d}\log f:= -V\frac{\partial \log f}{\partial\tau}+\frac{\partial
\log f}{\partial x}{\bf i} +\frac{\partial\log f}{\partial y}{\bf j}
+\frac{\partial\log f}{\partial z}{\bf k}\]
(notice that actually $\frac{\partial \log f}{\partial\tau}=0$).
Over $U_V\subset M_V$ we have a gauge induced by the above
orthonormal frame on the negative spinor bundle $\Sigma^-\vert_{U_V}$. 
In this gauge consider 't Hooft-like SU$(2)$ connections
$\nabla^-_{\lambda_0,\dots,\lambda_s}
\vert_{U_V}:=\dd +A^-_{\lambda_0,\dots,\lambda_s, U_V}$ on 
$\Sigma^-\vert_{U_V}$ with
\[A^-_{\lambda_0,\dots,\lambda_s, U_V}:={\rm Im}\frac{({\bf d} \log
f)\:\boldsymbol{\xi}}{ 2\sqrt{V}}.\] It was demonstrated in
\cite{ete-hau3} that these connections, parameterized by
$\lambda_0,\lambda_1,\dots,\lambda_s$ up to an overall scaling,
extend over $M_V$ and provide smooth, rapidly decaying anti-self-dual 
connections on $\Sigma^-$. They are irreducible if $\lambda_0>0$, 
are non-gauge equivalent and have trivial holonomy at infinity; hence satisfy 
the weak holonomy condition (\ref{holonomia}) and in particular the 
corresponding Chern--Simons invariants vanish. Consequently their energies are 
always integers equal to $e=n$ where $0\leq n\leq s$ is the number of 
non-zero $\lambda_i$'s with $i=1,\dots,s$. Therefore moduli spaces 
consisting of anti-instantons of these energies cannot be empty. $\Diamond$
\vspace{0.1in}

\begin{remark}\rm
It is reasonable to expect that the higher energy moduli spaces are also 
not empty. Furthermore, it is not clear whether fractionally charged, 
rapidly decaying irreducible instantons actually exist over 
multi-Taub--NUT or over more general ALF spaces at all.

Consider the family $\nabla^-_i$ defined by the
function $f_i(x)=\frac{\lambda_i}{\vert x-p_i\vert}$.
These are unital energy solutions which are reducible to U$(1)$ 
(in fact these are the only reducible points in
the family, cf. \cite{ete-hau3}). These provide as many as $s=1+b^2(M_V)$
non-equivalent reducible, rapidly decaying anti-self-dual solutions and if
$\nabla^-_i\vert_{U_V} =\dd + A^-_{i, U_V}$ then the connection
$\nabla^-\vert_{U_V}:=\dd +\sum\limits_i A^-_{i, U_V}$ decays rapidly,
is reducible and non-topological; hence it admits arbitrary rescalings 
yielding the strange continuous energy solutions mentioned in Section 2. 
Of course these generically rescaled connections violate the weak 
holonomy condition.
\end{remark}

\vspace{0.3in}

\centerline{\sc The Riemannian (or Euclidean) Schwarzschild space} 
\vspace{0.2in}

\noindent The underlying space is $M=S^2\times\R^2$. We have a 
particularly nice form of the
metric $g$ on a dense open subset $(\R^2\setminus \{ 0\})\times S^2\subset M$
of the Riemannian Schwarzschild manifold. It is convenient to use polar
coordinates $(r,\tau)$ on $\R^2\setminus\{0\}$ in the range $r\in
(2m,\infty )$ and $\tau\in [0,8\pi m)$, where $m>0$ is a fixed constant.
The metric then takes the form
\[\dd s^2=\left( 1-\frac{2m}{r}\right) \dd\tau^2+
\left( 1-\frac{2m}{r}\right)^{-1}\dd r^2+r^2\dd\Omega^2,\]
where $\dd\Omega^2$ is the line element of the round sphere.
In spherical coordinates $\theta\in (0,\pi)$ and $\varphi\in [0,2\pi)$ it 
is $\dd\Omega^2=\dd\theta^2+\sin^2\theta\:\dd\varphi^2$ on the open 
coordinate chart $(S^2\setminus(\{S\}\cup\{ N \} ))\subset S^2$.
Consequently the above metric takes the following form on the open, dense
coordinate chart $U:=(\R^2\setminus\{0\})\times
S^2\setminus(\{S\}\cup\{N\})\subset M$:
\[\dd s^2=\left( 1-\frac{2m}{r}\right)\dd\tau^2+
\left( 1-\frac{2m}{r}\right)^{-1}\dd r^2+r^2 (\dd\theta^2+\sin^2\theta
\dd\varphi^2).\]
The metric can be extended analytically to the whole $M$ as a complete 
Ricci flat metric however this time $W^\pm\not=0$. Nevertheless we obtain

\begin{theorem} Let $(M,g)$ be the Riemannian Schwarzschild manifold with
a fixed orientation. Let $\nabla_A$ be a smooth, rapidly decaying, 
self-dual SU$(2)$ connection on a fixed rank 2 complex vector bundle $E$ 
over $M$. Then $\nabla_A$ satisfies the weak holonomy condition and has 
non-negative integer energy. Let ${\cal M}(e,\Gamma )$ denote the framed 
moduli space of these connections which are moreover 
irreducible as in Theorem \ref{moduluster}. Then it is either empty or a 
manifold of dimension
\[{\dim\cal M}(e,\Gamma )=8e-3.\]
The moduli space with $e=1$ is surely non-empty.
\label{schwarz}
\end{theorem}

\noindent{\it Proof.} The metric is Ricci flat moreover a direct 
calculation shows 
that both $W^\pm$ satisfy the decay (\ref{gorbulet}) hence Theorem 
\ref{moduluster} applies in this situation as well. Furthermore, the 
asymptotical topology and the character variety of the space is
again $W\cong S^2\times S^1\times\R^+$ and $\chi (S^2\times S^1)\cong
[0,1)$ consequently the energy of a rapidly decaying instanton is integer 
as in Theorem \ref{kaloron}. We can again set the gauge in which all 
Chern--Simons invariants vanish and find moreover $X\cong S^2\times S^2$ 
yielding $b^-(X)=b^+(X)=1$; substituting these data into the dimension 
formula we get the desired result.

Regarding non-emptiness, very little is known. The apparently
different non-Abelian solutions found by Charap and Duff (cf. the
1-parameter family (I) in \cite{cha-duf}) are in fact all gauge equivalent 
\cite{tek} and provide a single rapidly decaying self-dual connection 
which is the positive chirality spin connection. It looks like 
$\nabla^+\vert_U:=\dd +A^+_U$ on $\Sigma^+\vert_U$ with
\[A^+_U:= \frac{1}{2}\sqrt{1-\frac{2m}{r}}\dd\theta{\bf i}
+\frac{1}{2}\sqrt{1-\frac{2m}{r}}\sin\theta\:\dd\varphi{\bf j}+
\frac{1}{2}\left(\cos\theta \dd\varphi
-\frac{m}{r^2}\dd\tau\right){\bf k}.\]
One can show that this connection extends smoothly over $\Sigma^+$ as
an ${\rm SO}(3)\times{\rm U}(1)$ invariant, irreducible, self-dual
connection of unit energy, centered around the 2-sphere in the origin.
$\Diamond$
\vspace{0.1in}

\begin{remark}\rm Due to its resistance against deformations
over three decades, it has been conjectured that this positive
chirality spin connection is the only unit energy instanton over the
Schwarzschild space (cf. e.g. \cite{tek}). However we can see now that in 
fact it admits a 2 parameter deformation. Would be interesting to find 
these solutions explicitly as well as construct higher energy irreducible 
solutions.

In their paper Charap and Duff exhibit another
family of ${\rm SO}(3)\times{\rm U}(1)$ invariant instantons
$\nabla_n$ of energies $2n^2$ with $n\in\Z$ (cf. solutions of type (II) in
\cite{cha-duf}). However it was pointed out in \cite{ete-hau1} that these 
solutions are in fact reducible to U$(1)$ and locally look like 
$\nabla_n\vert_{U^\pm}=\dd +A_{n,U^\pm}$ with
\[A_{n, U^\pm} :=\frac{n}{2}\left( (\mp 1+\cos\theta
)\dd\varphi -\frac{1}{r}\dd\tau\right){\bf k}\]
over the charts $U^\pm$ defined by removing the north or the south
poles from $S^2$ respectively. They extend
smoothly as {\it slowly decaying} reducible, self-dual connections 
over the bundles $L_n\oplus L^{-1}_n$ where $L_n$ is a line bundle 
with $c_1(L)=n$. Hence they are topological in contrast to the above 
mentioned Abelian instantons over the multi-Taub--NUT space. This 
constrains them to have discrete energy spectrum despite their slow 
decay. It is known that these are the only reducible SU$(2)$ instantons over 
the Riemannian Schwarzschild space \cite{ete-hau1,hau-hun-maz}.
\end{remark}

%%%%%%%%%%%%%%%%%%%%%%%%%%%%%%%%%%%%%%%%%%%%%%%%%%
%%%%%%%%%%%%%%%%%%%%%%%%%%%%%%%%%%%%%%%%%%%%%%%%%%

\section{Appendix}\label{secfive}

In this Appendix we shall prove Theorem \ref{index}; this proof is taken
from \cite{jar1}, and follows closely the arguments of \cite{gro-law}.
In the course of the proof we shall use the notation introduced in the
bulk of the paper.

The right hand side of the index formula in Theorem \ref{index} is called
the {\it relative topological index} of the operators $D_0$ and $D_1$:
\[{\rm Index}_t(D_1,D_0):= {\rm Index}_t\tilde{D}_1-{\rm
Index}_t\tilde{D}_0.\]
Notice that it can be computed in terms of the topology of the topological
extensions of the bundles $F_j$ and $F_j'$ to $X$, using the
Atiyah--Singer index theorem. Furthermore, as we will see below, this
quantity does not depend on how the operators $D_0$ and $D_1$ are extended
to $\tilde{D}_0$ and $\tilde{D}_1$ (see Lemma \ref{lem} below).

The first step in the proof of Theorem \ref{index} is the construction of
a new Fredholm operator $D^\prime_1$ as follows. Let $\beta_1$ and
$\beta_2$ be cut-off functions, respectively supported over $K$ and
$W=M\setminus K$, and define
\[D^\prime_1=\beta_1 D_1 \beta_1 + \beta_0 D_0 \beta_0.\]
Now, it is clear that $D^\prime_1|_W$ coincides with
$D_0\vert_W$. Furthermore, since
$\Vert (D^\prime_1-D_1)\vert_W\Vert <\kappa$ with $\kappa$
arbitrarily small, we know that ${\rm Index}_aD_1^\prime={\rm
Index}_aD_1$.
Our strategy is to establish the index formula for the pair $D_1'$ and
$D_0$. In order to simplify notation however, we will continue to denote
by $D_1$ and $D_0$ a pair of elliptic Fredholm operators which {\it coincide}
at infinity.

Now recall that if $D$ is any Fredholm operator over $M$, there is a
bounded, elliptic pseudo-differential operator $Q$, called the
{\it parametrix} of $D$, such that $DQ=I-S$ and $QD=I-S^\prime$, where $S$
and $S^\prime$ are compact {\it smoothing operators}, and $I$ is the 
identity operator. Note that neither $Q$ nor $S$ and $S^\prime$ are 
unique. In particular, there is a bounded operator $G$, called the {\it 
Green's operator} for $D$, satisfying $DG=I-H$ and $GD=I-H^\prime$, where 
$H$ and $H^\prime$ are finite rank {\it projection operators}; the image 
of $H$ is ${\rm Ker}D$ and the image of $H^\prime$ is ${\rm Coker}D$.

Let $K^H(x,y)$ be the Schwartzian kernel of the operator $H$. Its local
trace function is defined by ${\rm tr}[H](x)=K^H(x,x)$; moreover, these
are $C^\infty$ functions \cite{gro-law}. If $D$ is Fredholm, its (analytical)
index is given by
\begin{equation}
\label{ind}
{\rm Index}_aD=\dim{\rm Ker}D - \dim{\rm Coker}D =
\int\limits_M\left({\rm tr}[H]-{\rm tr}[H^\prime]\right)
\end{equation}
as it is well-known; recall that compact operators have smooth, square
integrable kernels. Furthermore, if $M$ is a closed manifold, we have
\cite{gro-law}
\[{\rm Index}_tD=\int\limits_M\left({\rm tr}[S]-{\rm tr}[S^\prime]\right)
.\]

Let us now return to the situation set up above. Consider the parametrices
and Green's operators ($j=0,1$)
\begin{equation}\label{ops}
\begin{array}{ccc}
\left\{ \begin{array}{l} D_jQ_j=I-S_j \\ Q_jD_j=I-S^\prime_j \end{array}
\right. &\ \ &
\left\{ \begin{array}{l} D_jG_j=I-H_j \\ G_jD_j=I-H^\prime_j .\end{array}
\right.
\end{array}
\end{equation}
The strategy of proof is to express both sides of the index formula of
Theorem \ref{index} in terms of integrals, as in (\ref{ind}); for its left
hand side, the {\it relative analytical index}, we have
\begin{equation}\label{inda2}
{\rm Index}_a(D_1,D_0) := {\rm Index}_aD_1 - {\rm Index}_aD_0 =
\int\limits_M\left({\rm tr}[H_1]-{\rm tr}[H_1^\prime]\right)-
\int\limits_M\left({\rm tr}[H_0]-{\rm tr}[H_0^\prime]\right) .
\end{equation}
For the relative topological index, we have the following

\begin{lemma}
\label{lem}
Under the hypothesis of Theorem \ref{index}, we have that
\begin{equation}\label{mid}
{\rm Index}_t(D_1,D_0) =
\int\limits_M\left({\rm tr}[S_1]-{\rm tr}[S_1^\prime]\right)-
\int\limits_M\left({\rm tr}[S_0]-{\rm tr}[S_0^\prime]\right) .
\end{equation}
\end{lemma}

\noindent{\it Proof.} Denote by $\widetilde{W}\subset X$ the 
compactification of
the end $W\subset M$. Extend $D_j$ ($j=0,1$) to operators $\tilde{D}_j$,
both defined over the whole $X$. Then the parametrices $\tilde{Q}_j$ of
$\tilde{D}_j$ are extentions of the parametrices $Q_j$ of $D_j$, and the
corresponding compact smoothing operators $\tilde{S}_j$ and $\tilde{S}_j'$
are extentions of $S_j$ and $S_j'$.

As explained above, we can assume that the operators $D_1$ and $D_0$
coincide at infinity. This means that $D_1|_W \simeq D_0|_W$, hence also
$\tilde{D}_1|_{\widetilde{W}} \simeq \tilde{D}_0|_{\widetilde{W}}$, and
therefore
\[S_1|_W \simeq S_0|_W ~~{\rm and} ~~ S_1'|_W \simeq S_0'|_W;\]
\[\tilde{S}_1|_{\widetilde{W}} \simeq \tilde{S}_0|_{\widetilde{W}} ~~{\rm
and} ~~ \tilde{S}_1'|_{\widetilde{W}} \simeq \tilde{S}_0'|_{\widetilde{W}}.\]
It follows that the operators $\tilde{S}_1 - \tilde{S}_0$ and
$\tilde{S}_1' - \tilde{S}_0'$ are supported on $K=M\setminus
W=X\setminus\widetilde{W}$ furthermore,
\[\tilde{S}_1 - \tilde{S}_0 = (S_1 - S_0)|_K ~~{\rm and}~~
\tilde{S}_1' - \tilde{S}_0' = (S_1' - S_0')|_K.\]
It follows that
\[{\rm Index}_t(D_1,D_0) = {\rm Index}_t\tilde{D}_1-{\rm
Index}_t\tilde{D}_0 =\]
\[=\int\limits_X\left({\rm tr}[\tilde{S}_1]-{\rm tr}[\tilde{S}_1']\right)
-\int\limits_X\left({\rm tr}[\tilde{S}_0]-{\rm tr}[\tilde{S}_0']\right) =
\int\limits_X\left( {\rm tr}[\tilde{S}_1]-{\rm tr}[\tilde{S}_0]\right)
-\int\limits_X\left( {\rm tr}[\tilde{S}_1']-{\rm tr}[\tilde{S}_0']\right)
=\]
\[=\int\limits_M\left( {\rm tr}[S_1]-{\rm tr}[S_0]\right)
-\int\limits_M\left( {\rm tr}[S_1']-{\rm tr}[S_0']\right) =
\int\limits_M\left({\rm tr}[S_1]-{\rm tr}[S_1^\prime]\right)-
\int\limits_M\left({\rm tr}[S_0]-{\rm tr}[S_0^\prime]\right)\]
as desired. $\Diamond$
\vspace{0.1in}

\noindent As we noted before, the proof of the lemma shows also that the
definition of the relative topological index is independent of the choice 
of extensions $\widetilde{D}_0$ and $\widetilde{D}_1$.

Before we step into the proof of Theorem \ref{index} itself, we must
introduce some further notation. Let $f:[0,1]\rightarrow [0,1]$ be a
smooth function such that $f=1$ on $[0,\frac{1}{3}]$, $f=0$ on
$[\frac{2}{3},1]$ and $f'\approx-1$ on $[\frac{1}{3},\frac{2}{3}]$.
Pick up a point $x_0\in M$
and let $d(x)={\rm dist}(x,x_0)$. For each $m\in\Z^*$, consider the
functions
\[f_m(x)=f\left(\frac{1}{m}\:{\rm e}^{-d(x)}\right) .\]
Note that ${\rm supp}\:\dd f_m^{\frac{1}{2}}\subset B_{\log
(\frac{3}{4m})}-B_{\log (\frac{3}{2m})}$ and
\begin{equation}
\Vert\dd f_m\Vert_{L^2}\leq\frac{c_9}{m}
\label{beta2}
\end{equation}
where $c_9=\left(\int_X{\rm e}^{-d(x)}\right)^{\frac{1}{2}}$. Here,
$B_r=\{x\in M\ |\ d(x)\leq r\}$, which is compact by the completeness of
$M$.
\vspace{0.1in}

\noindent {\it Proof of Theorem \ref{index}}. All we have to do is to show
that the right  hand sides of (\ref{inda2}) and (\ref{mid}) are equal.

In fact, let $U\subset V$ be small neighbourhoods of the diagonal
within $M\times M$ and choose $\psi\in C^\infty(M\times M)$
supported on $V$ and such that $\psi=1$ on $U$. Let $Q_j$ be
the operator whose Schwartzian kernel is
$K^{Q_j}(x,y)=\psi(x,y)K^{G_j}(x,y)$, where $G_j$ is the Green's
operator for $D_j$. Then $Q_j$ is a parametrix for $D_j$ for which
the corresponding smoothing operators $S_j$ and $S_j^\prime$, as in
(\ref{ops}), satisfy
\begin{equation}\label{traces} \begin{array}{ccc}
{\rm tr}[S_j]={\rm tr}[H_j] & {\rm and} & {\rm tr}[S_j^\prime]=
{\rm tr}[H_j^\prime]
\end{array} \end{equation}
where $H_j$ and $H_j^\prime$ are the finite rank projection operators
associated with the Green's operator $G_j$, as in (\ref{ops}).

But is not necessarily the case that the two parametrices $Q_0$ and
$Q_1$ so obtained must coincide at $W$. In order to fix that, we will glue
them with the common parametrix of $D_0|_W$ and $D_1|_W$, denoted $Q$ 
(with corresponding smoothing operators $S$ and $S'$),
using the cut-off functions $f_m$ defined above (assume that the base
points are contained in the compact set $K$). More precisely, for a 
section $s$
\[Q_j^{(m)}(s) = f_m^{\frac{1}{2}}\:Q_j(f_m^{\frac{1}{2}}\:s)+
(1-f_m)^{\frac{1}{2}}Q((1-f_m)^{\frac{1}{2}}s);\]
clearly, for each $m$, the operators $Q_0^{(m)}$ and $Q_1^{(m)}$ coincide 
at $W$. For the respective smoothing operators, we get (see 
\cite[Proposition 1.24]{gro-law})

\[\left\{ \begin{array}{ll}
S_j^{(m)}(s)= & f_m^{\frac{1}{2}}\:S_j(f_m^{\frac{1}{2}}\:s)+
(1-f_m)^{\frac{1}{2}}S((1-f_m)^{\frac{1}{2}}s)+\left( 
Q_j(f_m^{\frac{1}{2}}\:s)-
             Q((1-f_m)^{\frac{1}{2}}s\right)\dd f_m^{\frac{1}{2}},\\
              
S_j^{(m)\prime}(s)= & f_m^{\frac{1}{2}}\:S_j^\prime(f_m^{\frac{1}{2}}\:s)+
(1-f_m)^{\frac{1}{2}}S^\prime((1-f_m)^{\frac{1}{2}}s).
           \end{array}
\right.\]
Therefore
\[{\rm tr}[S_j^{(m)}]-{\rm tr}[S_j^{(m)'}]=
f_m^{\frac{1}{2}}\left({\rm tr}[S_j]-{\rm tr}[S'_j]\right) +
(1-f_m)^{\frac{1}{2}}\left({\rm tr}[S]-{\rm tr}[S']\right) +{\rm
tr}[(Q_j-Q)\dd f_m^{\frac{1}{2}}]\]
and
\[{\rm tr}[S_1^{(m)}]-{\rm tr}[S_1^{(m)'}]-
{\rm tr}[S_0^{(m)}]+{\rm tr}[S_0^{(m)'}]=\]
\[= f_m^{\frac{1}{2}}\left({\rm tr}[S_1]-{\rm tr}[S'_1]-{\rm tr}[S_0]+ 
{\rm tr}[S'_0]\right) +
{\rm tr}[(Q_1-Q)\dd f_m^{\frac{1}{2}}] - {\rm tr}[(Q_0-Q)\dd 
f_m^{\frac{1}{2}}]\]
so finally we obtain
\[{\rm tr}[S_1^{(m)}]-{\rm tr}[S_1^{(m)'}]-
{\rm tr}[S_0^{(m)}]+{\rm tr}[S_0^{(m)'}]=\]
\begin{equation}
\label{ff}
=f_m^{\frac{1}{2}}\left({\rm tr}[S_1]-{\rm tr}[S'_1]-
{\rm tr}[S_0]+ {\rm tr}[S'_0]\right) + {\rm tr}[(Q_1-Q_0)\dd 
f_m^{\frac{1}{2}}].
\end{equation}

We must now integrate both sides of (\ref{ff}) and take limits as
$m\rightarrow\infty$. For $m$ sufficiently large, ${\rm
supp}(1-f_m)\subset W$, hence the left hand side of
the identity (\ref{ff}) equals the relative topological index
${\rm Index}_t(D_1,D_0)$,
by Lemma \ref{lem}. On the other hand, the first summand inside on the
right hand side of (\ref{ff})
equals ${\rm Index}_a(D_1,D_0)$ by (\ref{traces}) and (\ref{inda2}).
Thus, it is enough to show that the integral of the last two terms on the
right hand side of (\ref{ff}) vanishes as $m\rightarrow\infty$. Indeed,
note that
\[{\rm tr}[(Q_1-Q_0)\dd f_m^{\frac{1}{2}}]=
\dd f_m^{\frac{1}{2}}\:{\rm tr}[Q_1-Q_0]\]
hence, since ${\rm supp}(\dd f_m)\subset W$ for sufficiently large $m$
and using also (\ref{beta2}), it follows that
\begin{equation}\label{ineq}
\int\limits_W{\rm tr}[(Q_1-Q_0)\dd f_m^{\frac{1}{2}}]\leq
\frac{c_9}{m}\int\limits_W{\rm tr}[G_1-G_0]\rightarrow 0~ {\rm as}~
m\rightarrow\infty
\end{equation}
if the integral on the right hand side of the above inequality is finite.
Indeed, let $D=D_1|_W=D_0|_W$; from the parametrix equation, we have
$D((G_1-G_0)|_W)=(H_1-H_0)|_W.$
Observe that ${\cal H}={\rm Ker}((H_1-H_0)|_W)$ is a closed
subspace of finite codimension in $L^2_{j+1,\Gamma_1}(F_1|_W)$. 
Moreover ${\cal H}\subseteq{\rm Ker}D$; thus, $(G_1-G_0)|_W$ has finite
dimensional range and hence it is of trace class i.e., the integral
on the right hand side of inequality (\ref{ineq}) does converge
(see also \cite[Lemma 4.28]{gro-law}). This concludes the proof.
$\Diamond$
\vspace{0.1in}

\noindent{\bf Acknowledgement.}
The authors thank U. Bunke, S. Cherkis, T. Hausel and D. N\'ogr\'adi for 
careful reading of the manuscript and making several 
important comments. We also thank the organizers of the 2004 March 
AIM-ARCC workshop on $L^2$ harmonic forms for bringing the team together. 
G.E. is supported by the CNPq grant No. 150854/2005-6 (Brazil) and the OTKA 
grants No. T43242 and No. T046365 (Hungary). M.J. is partially supported by 
the CNPq grant No. 300991/2004-5.

\end{document}